\documentclass[11 pt]{report}

\usepackage{amsmath}
\usepackage{amsfonts}
\usepackage{latexsym}
\usepackage{graphicx}
\usepackage{ifpdf} 
\ifpdf 
\fi 
\usepackage[all]{xy}
\usepackage{rotating}
\usepackage{theorem}
\usepackage[mathcal]{euscript}
\usepackage{amscd}
\usepackage{amssymb}
\usepackage{comment}
\usepackage{eufrak}
\usepackage{epsfig}
\usepackage{subfigure}
\usepackage{float}
\usepackage{epic}
\usepackage{eepic}
\usepackage{fullpage}
\usepackage{mathdots}

\newcommand{\QED}{\hfill $\blacksquare$}

\renewcommand{\H}{\mathcal{H}}
\renewcommand{\L}{\mathcal{L}}

\newcommand{\A}{\mathcal{A}}

\newcommand{\ZZ}{\mathbb{Z}}

\newcommand{\ST}{\boldsymbol{s}}
\newcommand{\UU}{\boldsymbol{u}}
\newcommand{\Y}{\boldsymbol{y}}

\newtheorem{theorem}{Theorem}[section]
\newtheorem{definition}[theorem]{Definition}
\newtheorem{lemma}[theorem]{Lemma}

\newtheorem{prop}[theorem]{Proposition}

\newtheorem{example}[theorem]{Example}
\newtheorem{remark}[theorem]{Remark}

\title{Structural properties of acyclic heaps with applications to Kazhdan--Lusztig theory}
\author{\textbf{Brent  G. Pohlmann}\\
Ph.D., University of Colorado at Boulder, 2008}

\begin{document}
\maketitle
\setcounter{page}{1}
\pagenumbering{roman}
\begin{abstract}  
\begin{center}
\textbf{Structural properties of acyclic heaps of pieces with Kazhdan--Lusztig theory}

Thesis directed by Prof. Richard M. Green
\end{center}

Let  {\it W} be an arbitrary Coxeter group with generating set {\it S} of involutions. A reduced expression for an element $w \in W$ is a minimal length word in {\it S} that represents {\it w}. The set $W_c$ of fully commutative elements are characterized by the property that any reduced expression for $w \in W_c$ can be obtained from any other via iterated commutations of adjacent generators. If $w \in W_c$ and $sw \notin W_c$ for some $s \in S$, then we say $sw$ is weakly complex.

Star reducible Coxeter groups are a class of Coxeter groups whose fully commutative elements have a particularly nice property. Star reducible Coxeter groups contain the finite Coxeter groups as a subclass.

A heap is an isomorphism class of labelled posets. Each heap is equipped with a set of edges and a set of vertices. Green defined a linear map $\partial_E$ which sends each edge of a heap {\it E} to a linear combination of vertices. If $v \in \textrm{Im}\partial_E$, we call {\it v} a boundary vertex. If $e_0$ is an edge of {\it E} and $\partial_E(e_0) = v$ for a vertex {\it v}, then we call {\it v} an effective boundary vertex. If $\partial_E(e_0) = v_1 + v_2$, then $v_1$ and $v_2$ are said to be linearly equivalent. A result by Stembridge states that every fully commutative element has a unique heap.

The main result we will prove here is that in the heap of a fully commutative element in a star reducible Coxeter group, every boundary vertex is linearly equivalent to an effective boundary vertex. 

We use the main result to prove another theorem concerning the basis elements  $\widetilde{t}_w : w \in W_c$ of a generalized Temperley--Lieb algebra. The theorem we prove 
allows an inductive computation of the $c$-basis which is a linear combination of the $\widetilde{t}_w : w \in W_c$. Furthermore, the $c$-basis can be shown to have nonnegative structure constants, that is, structure constants that are Laurent polynomials with nonnegative coefficients. One of the reasons this is interesting is that in many cases, these structure constants are also structure constants for the well known Kazhdan--Lusztig basis $C'_w$, whose positivity is generally very difficult to prove.
\end{abstract}

\tableofcontents
\newpage
\pagenumbering{arabic}
\setcounter{page}{1}

\chapter{Introduction}
A Coxeter group {\it W} is a group with finite generating set {\it S} and presentation given by
\begin{displaymath}
W = \langle S \  |  \ (st)^{m(s,t)} = 1 \ \textnormal{for} \ m(s,t) < \infty \rangle.
\end{displaymath}
Each $w \neq 1$ in {\it W} can be written in the form $w = s_1s_2 \cdots s_r$ for some $s_i \in S$. If {\it r} is as small as possible we call {\it r} the length of {\it w}, written $\ell(w)$.
A reduced expression for an element $w \in W$ is a minimal length word in {\it S} that represents {\it w}. The set $W_c$ of fully commutative elements of {\it W} is characterized by the property that any reduced word for $w \in W_c$ can be obtained from any other via iterated applications of short braid relations, that is, relations of the form $st = ts$, where $s, t \in S$. For example, if {\it w} is a product of commuting generators from {\it S}, then {\it w} is fully commutative. If $w \in W_c$ and $sw \notin W_c$ for some $s \in S$, then we say {\it sw} is weakly complex.
Green \cite{Green:P} defined the star reducible Coxeter groups to be those Coxeter groups for which every fully commutative element is equivalent to a product of commuting generators by a sequence of length-decreasing star operations. The star reducible Coxeter groups are defined and listed in Section 2.1. 

Denote by $\mathcal{H}_q = \mathcal{H}_q(X)$ the Hecke algebra associated to {\it W}. This is a $\mathbb{Z}[q, q^{-1}]$-algebra with a basis consisting of (invertible) elements $\{T_w : w \in W \}$,  and an associative multiplication. For our purposes, we extend the scalars of $\mathcal{H}_q$ by setting $v^2 = q$ and writing $\mathcal{H} = \mathcal{A}\otimes_{\mathbb{Z}[q, q^{-1}]}\H_q$ where $\mathcal{A} = \ZZ[v, v^{-1}]$. We write $\mathcal{A^+}$ and $\mathcal{A^-}$ for $\ZZ[v]$ and $\ZZ[v^{-1}]$, respectively. We also define a scaled version of the {\it T}-basis, $\{\widetilde{T}_w : w \in W\}$, where $\widetilde{T}_w := v^{-\ell(w)}T_w$. In \cite{Kazhdan;Lusztig:A}, Kazhdan and Lusztig defined the bases \{$C_w : w \in W$\} and  \{$C'_w : w \in W$\} for $\mathcal{H}$. These Kazhdan--Lusztig bases are constructed from the bases \{$T_w : w \in W$\}. An equation relating the $C'_w$-basis and the $\widetilde{T}_w$-basis is given by
\begin{displaymath}
C'_w = \widetilde{T}_w + \sum_{\substack{y < w}}A_{y,w} \widetilde{T}_y
\end{displaymath}
where $<$ is the Bruhat order on $W$ and $A_{y,w} \in v^{-1}\mathcal{A}^-$. Following \cite[Section 11.1]{M.Geck;G.Pfeiffer:A}, we denote the coefficients of $\widetilde{T}_y$ in $C'_w$ by $P^*_{y,w}(q)$. The Kazhdan--Lusztig polynomial $P_{y,w}$ is then given by $v^{\ell(w) - \ell(y)}P^*_{y,w}(q)$. The coefficient of the $v^{-1}$ term in $P^*_{y,w}(q)$ is $\mu(y, w)$, which is very difficult to compute efficiently, even for moderately small groups.

Let $J(X)$ be the two-sided ideal of $\H$ generated by the elements 
\begin{displaymath}
\sum_{\substack{w \in \langle s, t \rangle}} T_w,
\end{displaymath}
where $(s,t)$ runs over all pairs of elements of {\it S} such that $2 < m(s,t) < \infty$ and $\langle s, t \rangle$ is the parabolic subgroup generated by {\it s} and {\it t}.
Following Graham \cite{Graham:A}, we define the generalized Temperley--Lieb algebra $TL(X)$ to be the quotient $\mathcal{A}$-algebra $\H(X) / J(X)$ and denote the corresponding epimorphism of algebras by $\theta : \H(X) \rightarrow TL(X)$. For $\widetilde{T}_w \in \H$, we have $\theta(\widetilde{T}_w) = \widetilde{t}_w$ and define the $\mathcal{A^-}$-submodule $\L$ of $TL(X)$ to be that generated by the $\{\widetilde{t}_w : w \in W_c\}$. 

Green and Losonczy \cite{Green;Losonczy:A} proved that for each $w \in W_c$ there exists a unique $c_w \in TL(X)$ such that
\begin{displaymath}
c_w = \widetilde{t}_w + \sum_{\substack{y < w \\ y \in W_c}} a_{y,w}\widetilde{t}_y,
\end{displaymath}
where $<$ is the Bruhat order on {\it W}, and $a_{y,w} \in \mathcal{A^-}$ for all {\it y}. We note here that the $c_w \in TL(X)$ are analogous to $C'_w \in \H$ and the direct connection, which we explore later, is not obvious. 
 
A heap is an isomorphism class of labelled posets satisfying certain axioms. Each heap is equipped with a set of edges and a set of vertices. In \cite{Green:B}, Green defined a linear map $\partial_E$ which sends each edge of {\it E} to a linear combination of vertices. We say {\it v} is a  boundary vertex if $v \in \textrm{Im} \ \partial_E$ and that {\it v} is an effective boundary vertex if $\partial_E(e_1) = v$ for an edge $e_1$ from {\it E}. Two boundary vertices $v_1$ and $v_2$ are said to be linearly equivalent if $\partial_E(e_1) = v_1 + v_2$ for some edge $e_1$ from {\it E}.
In \S{3.4} we prove that in the heap of a fully commutative element in a star reducible Coxeter group, every boundary vertex is linearly equivalent to an effective boundary vertex. This is a subtle structural property of the heap and the proof is combinatorial in nature. 

As an application of Theorem 3.4.1, we prove that if $x \in W$ is weakly complex for {\it W}  star reducible, we have $\widetilde{t}_x \in v^{-1}\mathcal{L}.$ This property of the $\widetilde{t}_x$, which Green calls Property W  in \cite{Green:K}, allows the inductive computation of the $c$-basis using the formula
\begin{displaymath}
c_sc_w = \left \{ \begin{array}{ll}
(v + v^{-1})c_w, & \textrm{if $\ell(sw) < \ell(w)$},\\
c_{sw} + \sum_{\substack{sy < y}}\mu(y, w)c_y & \textrm{if $\ell(sw) > \ell(w)$}.
\end{array} \right.
\end{displaymath}
Under Property W, the $c$-basis can be shown to have nonnegative structure constants, that is, structure constants that are Laurent polynomials with nonnegative coefficients. One of the reasons this is interesting is that in many cases, these structure constants are also structure constants for the Kazhdan--Lusztig basis $C'_w$, whose positivity is generally very difficult to prove.

\chapter{Coxeter Group Theory}
\section{Introduction}
In \S{2} we review general Coxeter group theory. For a complete review of this material, the reader is referred to \cite{Humphreys:A} and \cite{Bjorner;Brenti:A}.

A Coxeter system $(W, S)$ consists of a group {\it W} with distinguished (finite) set of generating involutions {\it S} and presentation given by 
\begin{displaymath}
W = \langle S \  |  \ (st)^{m(s,t)} = 1 \ \textnormal{for} \ m(s,t) < \infty \rangle,
\end{displaymath}
where $m(s,s) = 1$. (It turns out that the $s \in S$  are distinct as group elements, and that $m(s,t)$ is the order of $st$.) Since the generators $s \in S$ have order 2 in {\it W}, each $w \neq 1$ in {\it W} can be written in the form $w = s_1s_2 \cdots s_r$ for some $s_i \in S$. If {\it r} is as small as possible we call {\it r} the {\it length} of {\it w}, written $\ell(w)$. A product $w_1w_2 \cdots  w_n$ of elements $w_i \in W$ is called {\it reduced} if $\ell(w_1w_2 \cdots  w_n) = \sum_{i} \ell(w_i)$. We reserve the terminology {\it reduced expression} for reduced products $w_1w_2 \cdots w_n$ in which every $w_i \in S$. To specify a Coxeter system $(W,S)$ we can specify a finite set $S$ and draw an undirected graph {\it X} with set $S$ as vertex set, joining vertices $s$ and $t$ by an edge labelled $m(s,t)$ whenever this number is at least 3. If distinct vertices $s$ and $t$ are not joined, it is understood that $m(s,t) =2$. As a simplifying convention, the label $m(s,t) = 3$ is usually omitted. We call {\it X} a {\it Coxeter graph}. The connected graph in Figure 1 is called a Coxeter graph of type $B_4$.

\begin{center}

\includegraphics{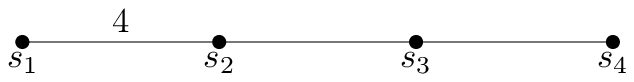}

FIG. 1. Coxeter graph of type $B_4$.
\end{center}

We call an element $ w \in W$ {\it complex} if it can be written as a reduced product $x_1w_{st}x_2$, where $x_1, x_2 \in W$ and $w_{st}$ is the longest element of some rank 2 parabolic subgroup $\langle s,t \rangle$ such that $s$ and $t$ correspond to adjacent nodes in the Coxeter graph. An element $w \in W$ is said to be $weakly \ complex$ if (a) it is complex and (b) it is of the form $w=su$, where $s \in S$ and $u$ is not complex. We write
\begin{displaymath}
\mathcal{L}(w) = \{s \in S : \ell(sw) < \ell(w)\}
\end{displaymath}
\textnormal{and}
\begin{displaymath}
\mathcal{R}(w) = \{s \in S : \ell(ws) < \ell(w)\}.
\end{displaymath}
The set $\mathcal{L}(w)$ (respectively, $\mathcal{R}(w)$) is called the {\it left} (respectively, {\it right}) {\it descent  set} of {\it w}.

\begin{example}
The reduced element $w = s_1s_2s_1s_2s_3$ in type $B_4$ is complex: the subexpression $w' = s_1s_2s_1s_2$ is the longest element of the rank 2 parabolic subgroup $\langle s_1, s_2  \rangle$. If we delete the first occurrence of $s_1$ in $w$ we are left with an element that is not complex. Thus {\it w} is also weakly complex, with $s = s_1$ in the definition above.
\end{example}

Denote by $W_c(X)$ the set of all elements of $W$ that are not complex. The elements of $W = W_c(X)$ are the {\it fully commutative} elements of \cite{Stembridge:B}; they are characterized by the property that any two of their reduced expressions may be obtained from each other by repeated commutation of adjacent generators.

Denote by $\mathcal{H}_q = \mathcal{H}_q(X)$ the Hecke algebra associated to $W$. This is a $\mathcal{A}$-algebra with a basis consisting of (invertible) elements $T_w : w \in W$, satisfying
\begin{displaymath}
T_sT_w = \left \{ \begin{array}{ll}
T_{sw}, & \textrm{if $\ell(sw) > \ell(w)$}\\
qT_{sw} + (q-1)T_w, & \textrm{if $\ell(sw) < \ell(w)$}
\end{array} \right.
\end{displaymath}
where $\ell$ is is the length function on the Coxeter group $W$, $w \in W$, and $s \in S$. Note that for $s \in S$, we have ${T}^{2}_s = (q - 1)T_s + qT_1.$

If $w = s_1s_2 \cdots s_r$ is a reduced expression for {\it w}, then $T_w = T_{s_1}T_{s_2} \cdots T_{s_r}$. Since $\mathcal{H}$ can also be shown to be associative, the product $T_xT_{x'}$ is defined for all $x, x' \in W.$

Using the above relations we can show that for all $s \in S$:
\begin{displaymath}
T^{-1}_s = q^{-1}T_s - (1 - q^{-1})T_1.
\end{displaymath}
Therefore every $T_w$ is invertible in $\H$.

For our purposes, we extend the scalars of $\mathcal{H}_q$ by setting $v^2 = q$ and define $\mathcal{H} = \mathcal{A}\otimes_{\mathbb{Z}[q, q^{-1}]}\H_q$. We also define a scaled version of the $T$-basis, $\{\widetilde{T}_w : w \in W\}$, where $\widetilde{T}_w := v^{-l(w)}T_w$. We will write $\A^+$ and $\A^-$ for $\mathbb{Z}[v]$ and $\mathbb{Z}[v^{-1}]$, respectively. We denote the $\mathbb{Z}$-linear ring homomorphism $\A \longrightarrow \A$ exchanging $v$ and $v^{-1}$ by $\bar{\ \ } $. We can extend $\bar{\ \ } $ to a ring automorphism of $\mathcal{H}$ by the condition that 
\begin{displaymath}
\overline{\sum_{w \in W} a_w {T_w}} = \sum_{w \in W} \overline{a_w} {T^{-1}_{w^{-1}}}
\end{displaymath}
where the $a_w$ are elements of $\mathcal{A}$.

Let $J(X)$ be the two-sided ideal of $\H$ generated by the elements 
\begin{displaymath}
\sum_{\substack{w \in \langle s, t \rangle}} T_w,
\end{displaymath}
where $(s,t)$ runs over all pairs of elements of {\it S} that correspond to adjacent nodes in the Coxeter graph, and $\langle s, t \rangle$ is the (finite) parabolic subgroup generated by {\it s} and {\it t}.
Following Graham \cite{Graham:A}, we define the {\it generalized Temperley--Lieb algebra} $TL(X)$ to be the quotient $\mathcal{A}$-algebra $\H(X) / J(X)$. We denote the corresponding epimorphism of algebras by $\theta : \H(X) \rightarrow TL(X)$. Let $t_w$ (respectively, $\widetilde{t}_w$) denote the image in $TL(X)$ of the basis element $T_w$ (respectively, $\widetilde{T}_w$) of $\H$. 

In \cite{Kazhdan;Lusztig:A}, Kazhdan and Lusztig defined the bases \{$C_w : w \in W$\} and  \{$C'_w : w \in W$\} for $\mathcal{H}$. These Kazhdan--Lusztig bases are constructed from the bases \{$T_w : w \in W$\}. The following theorem can be viewed as a restatement of \cite[1.1.c]{Kazhdan;Lusztig:A}. 
\begin{theorem} (Kazhdan, Lusztig)
For each $w \in W$, there exists a unique $C'_w \in \mathcal{H}$ such that both $\overline{C'_w} = C'_w$ and
\begin{displaymath}
C'_w = \widetilde{T}_w + \sum_{\substack{y < w}}A_{y,w} \widetilde{T}_y,
\end{displaymath}
where $<$ is the Bruhat order on $W$ and $A_{y,w} \in v^{-1}\mathcal{A}^-$ are certain polynomials in $v^{-1}$. 
\end{theorem}
The basis $\{C_w\}_{w \in W}$ is closely related to $\{C'_w\}_{w \in W}$, and it suffices to understand one of them. 

We define the $\mathcal{A}^-$ submodule $\mathcal{L}$ of $TL(X)$ to be that generated by $\{\widetilde{t}_w : w \in W_c\}$. We define $\pi: \L \longrightarrow \L/ v^{-1}\L$ to be the canonical $\mathbb{Z}-$linear projection. It was proved in \cite[Lemma 1.4]{Green;Losonczy:A} that the ideal $J(X)$ is fixed by $\bar{\ \ }$, so $\bar{\ \ }$ induces an involution on $TL(X)$ that sends $v$ to $v^{-1}$ and ${t}_w$ to $t_{w^{-1}}^{-1}$. We denote this map also by $\bar{\ \ }$.

The following theorem establishes a canonical basis for $TL(X)$ in terms of the $\widetilde{t}$-basis. It comes from \cite{Green;Losonczy:A}.

\begin{theorem}
For each $w \in W_c$, there exists a unique $c_w \in TL(X)$ such that both $\overline{c_w} = c_w$ and $\pi(c_w) = \pi(\widetilde{t_w})$. Furthermore, we have 
\begin{displaymath}
c_w = \widetilde{t}_w + \sum_{\substack{y < w \\ y \in W_c}} a_{y,w}\widetilde{t}_y,
\end{displaymath}
where $<$ is the Bruhat order on $W$, and $a_{y,w} \in v^{-1}\A^-$ for all $y$.
\end{theorem}

For later purposes, we define the following sublattices of the $\A^-$-lattice $\mathcal{L}$.
\begin{definition}
Let $s \in S$. We define $\mathcal{L}^s_L$ to be the free $\A^-$-module with basis
\begin{displaymath}
\{\widetilde{t}_w : w \in W_c, \ sw < w\} \cup \{v^{-1}\widetilde{t}_w : w \in W_c, \ sw > w\}.
\end{displaymath}
\end{definition}
Similarly, we define a free $\mathcal{A^-}$-module $\mathcal{L}^s_R$. 

Star operations are of key importance to this thesis. These were introduced in the simply laced case in \cite{Kazhdan;Lusztig:A} and in general in \cite{Lusztig:A}.

\begin{definition}
Let $W=W(X)$  be any Coxeter group and let $I = \{s, t \} \subseteq S$ be a pair of noncommuting generators whose product has order $m = m(s,t)$ (where $\infty$ is allowed). Let $W^I$ denote the set of all $w \in W$ satisfying $\mathcal{L}(w)  \cap I = \emptyset$. Standard properties of Coxeter groups \cite{Humphreys:A} show that any element $w \in W$ may be uniquely written as $ w = w_I w^I$, where $w_I \in W_I = \langle s, t \rangle$, $w^I \in W^I$ and $\ell(w) = \ell(w_I) + \ell(w^I)$. There are four possibilities for the elements $w \in W$:\\
(i) $w$ is the shortest element in the coset $W_Iw$, so $w_I =1$ and $w \in W^I$;\\
(ii) $w$ is the longest element in the coset $W_Iw$, so $w_I$ is the longest element of $W_I$ (which can only happen if $W_I$ is finite);\\
(iii) $w$ is one of the $(m-1)$ elements $ sw^I, tsw^I, stsw^I, \ldots$;\\
(iv) $w$ is one of the $(m-1)$ elements $ tw^I, stw^I, tstw^I, \ldots$.
\end{definition}

The sequences appearing in (iii) and (iv) are called {\it (left) \{s,t\}-strings}, or {\it strings} if the context is clear. If {\it x} and {\it y} are two elements of an $\{s, t\}$-string such that $\ell(x)= \ell(y)-1$, we call the pair $\{x,y\}$  {\it left \{s,t\}-adjacent}, and we say that {\it y} is {\it left star reducible} to {\it x}.
The above concepts all have right-handed counterparts, leading to the notion of {\it right \{s,t\}-adjacent} and {\it right star reducible} pairs of elements, and coset decompositions $(^Iw)(_Iw)$.
If there is a (possibly trivial) sequence
\begin{displaymath}
x = w_0, w_1, \ldots, w_k=y
\end{displaymath}
where, for each $0 \leq i <  k, w_{i+1}$ is left star reducible or right star reducible to ${w_i}$ with respect to some pair $\{s, t\}$, we say that $y$ is $star \ reducible \ to \ x$. Because star reducibility decreases length, this relation is antisymmetric and thus defines a partial order on {\it W}.
If {\it w} is an element of a left $\{s,t\}$-string, $S_w$, we have $\{\ell(sw), \ell(tw)\} = \{\ell(w) - 1, \ell(w) + 1\}$; let us assume without loss of generality that $sw$ is longer than $w$ and $tw$ is shorter. If $sw$ is an element of $S_w$, we define $^*w = sw$; if not, $^*w$ is undefined. If $tw$ is an element of $S_w$, we define $_*w= tw$; if not, $_*w$ is undefined.
There are also obvious right handed analogues to the above concepts, so the symbols $w^*$ and $w_*$ may be used with the analogous meanings.

\begin{example}
Let {\it W} the Coxeter group of type $B_4$ and let $I = \{s_1, s_2\}$, where $m(s_1, s_2) = 4$. If we consider the element $w = s_2 s_1$, we have
\begin{displaymath}
_*w = s_1, \ ^*w = s_1s_2s_1, \ w_* = s_2 \ \textnormal{and} \ w^* = s_2s_1s_2. 
\end{displaymath}
If $x = s_1s_2s_1$ then $^*x$ is undefined and $x^*$ is undefined; if $x = s_2$ then $_*x$ and $x_*$ are undefined.
\end{example}

\begin{definition}
We say that a Coxeter group $W(X)$, or its Coxeter graph {\it X}, is {\it star reducible} if every element of $W_c$ is star reducible to a product of commuting generators from {\it S}. Either $X$ is a complete graph $K_n$ with all labels $m(s, t) \geq 3$, or $X$ appears in the list below, which comes from \cite{Green:P}.
\end{definition}

\vspace{0.5 in}

\begin{center}

\includegraphics{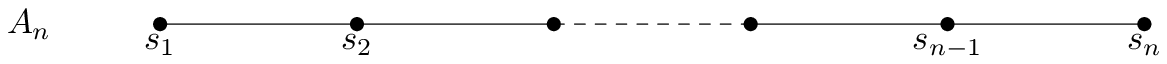}

\vspace{0.5 in}

\includegraphics{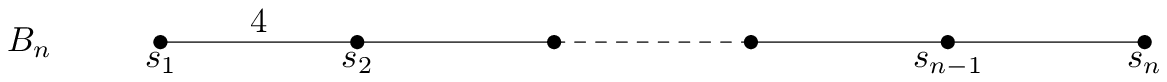}

\vspace{0.5 in}

\includegraphics{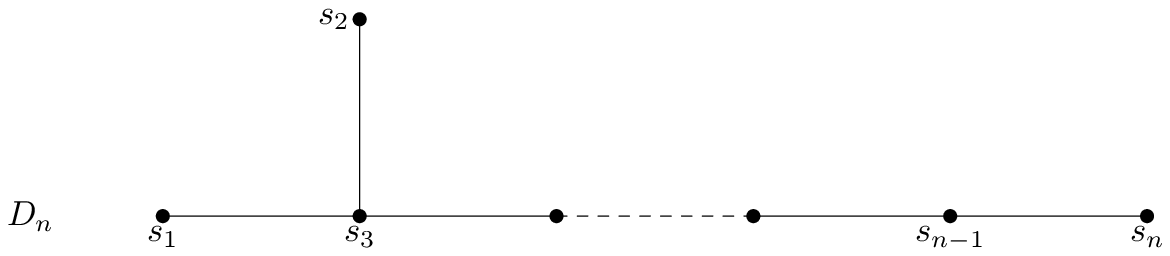}

\vspace{0.5 in}

\includegraphics{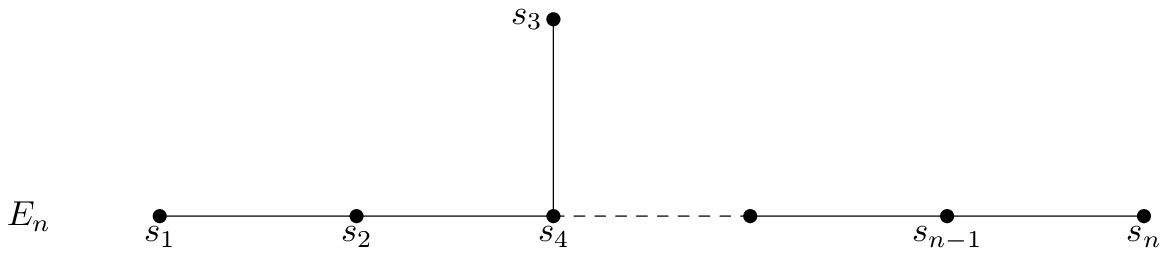}

\vspace{0.5 in}

\includegraphics{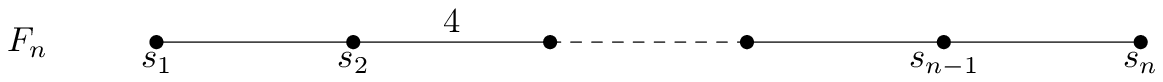}

\vspace{0.5 in}

\includegraphics{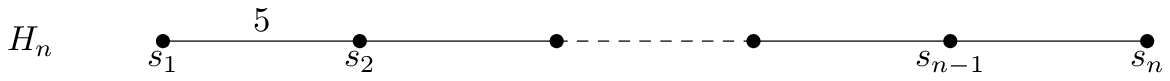}

\vspace{0.5 in}

\includegraphics{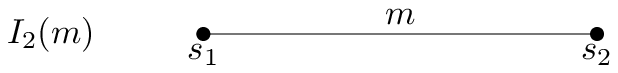}

\vspace{0.5 in}

\includegraphics{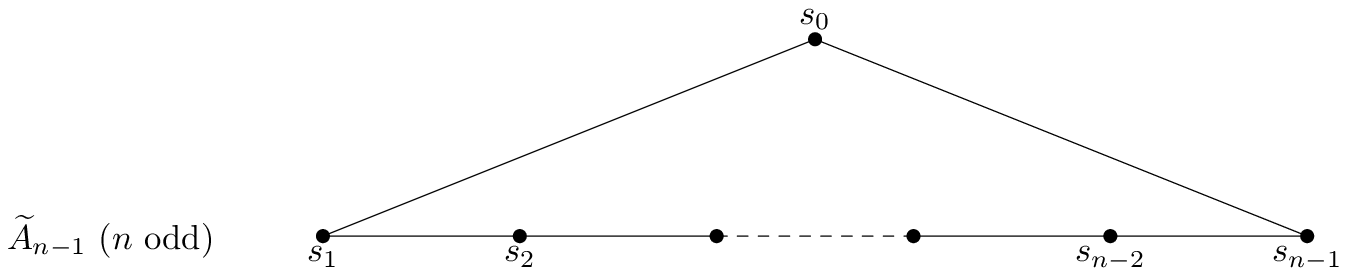}

\vspace{0.5 in}

\includegraphics{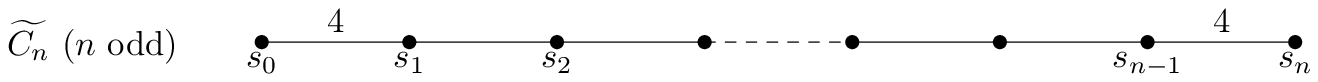}

\vspace{0.5 in}

\includegraphics{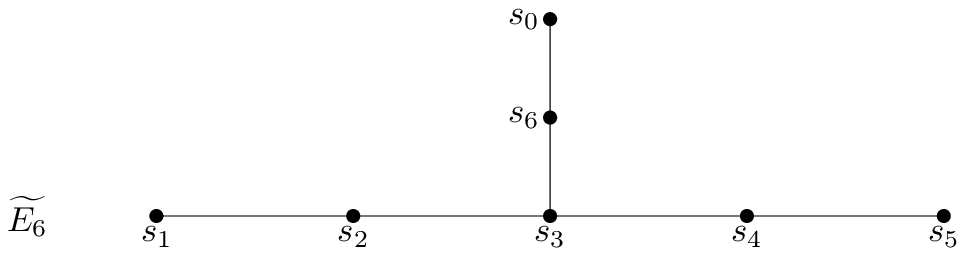}

\vspace{0.5 in}

\includegraphics{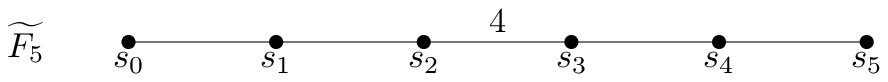}

\vspace{0.5 in} 

FIG. 2 The star reducible Coxeter groups
\end{center}

\begin{definition}
If we have $m(s,t) \leq 3$ for all $s, t \in S$, we say {\it W} is {\it simply laced}.
\end{definition}

The following proposition is a subset of \cite[Proposition 4.10]{Green:K}. Part (iii) of the original had a slight change and part (iv) was completely wrong. Note that a right-handed analogue of the proposition  also holds. 
\begin{prop}
Suppose that {\it s} and {\it t} are noncommuting generators of the Coxeter group {\it W}, and that $\widetilde{t}_x \in \mathcal{L}^u_L$ whenever {\it x} is weakly complex, $ux \in W_c$ and $u \in S$. Then we have :\\
(i)
\begin{displaymath} 
\widetilde{t}_s\widetilde{t}_w \in \left\{ \begin{array}{ll} 
v\mathcal{L}^s_L & \textrm{if $sw < w,$}\\ 
 \mathcal{L}^s_L & \textrm{if $sw > w;$}\\ 
\end{array} \right. 
\end{displaymath}
\end{prop}
(ii) $\widetilde{t}_s\mathcal{L} \cap  \mathcal{L} \subseteq  \mathcal{L}^s_L;$\\
(iii) $\widetilde{t}_s\mathcal{L}^t_L \subseteq  \mathcal{L}^s_L.$ \QED

\chapter{Heaps of Pieces}
\section{Introduction}
In Section 3.1, we introduce the basic properties of heaps. We will tend to follow Viennot's notation \cite{Viennot:A}.
\begin{definition} %3.1.1
Let {\it P} be a set equipped with a symmetric and reflexive binary relation $\mathcal{C}$. The elements of {\it P} are called {\it pieces}, and the relation $\mathcal{C}$ is called the {\it concurrency relation}.
A {\it labelled heap} with pieces in $P$ is a triple $(E, \leq, \varepsilon)$ where $(E, \leq)$ is a finite (possibly empty) partially ordered set with order relation denoted by $\leq$ and $\varepsilon$ is a map $\varepsilon : E \rightarrow P$ satisfying the following two axioms:
\end{definition}
 \noindent (i) for every $a, b \in E$ such that $\varepsilon{(a)} \ \mathcal{C} \ \varepsilon{(b)}$, {\it a} and {\it b} are comparable in the order $\leq$; \\
\noindent (ii) the order relation $\leq$ is the reflexive and transitive closure of the relation $\leq_\mathcal{C}$ such that for all $a, b \in E, a \leq_\mathcal{C} b$ if and only if both $a \leq b$ and $\varepsilon{(a)} \ \mathcal{C} \ \varepsilon{(b)}.$

\begin{definition}%3.1.2
Let  $(E, \leq, \varepsilon)$ and let $(E', \leq', \varepsilon')$ be two labelled heaps with pieces in $P$ and the same concurrency relation, $\mathcal{C}$. An isomorphism $\phi : E \rightarrow E'$ of posets is said to be an {\it isomorphism of labelled posets} if $\varepsilon = \varepsilon' \circ \phi$.
\end{definition}

\begin{definition}%3.1.3
A {\it heap} of pieces in {\it P} with concurrency relation $\mathcal{C}$ is an equivalence class of labelled heaps under labelled poset isomorphism. The set of such heaps is denoted $H(P, C)$. We denote the heap corresponding to the labelled heap $(E, \leq, \varepsilon)$ by $[E, \leq, \varepsilon]$. We will sometimes abuse notation by referring to a heap $[E, \leq, \varepsilon]$ simply as $E$.
\end{definition} See Example 3.1.6 for an example of a labelled heap.

\begin{definition}%3.1.4
The {\it concurrency graph} associated to the class of heaps $H(P, \mathcal{C})$ is the graph $X$ whose vertices are the elements of {\it P} and for which there is an edge from $v \in P$ to $w \in P$ if and only if $v \neq w$ and $v \ \mathcal{C} \ w$. We will write $H(X)$ to mean $H(P, C)$.
\end{definition}

We are particularly interested in heaps arising from fully commutative elements in Coxeter groups as studied by Stembridge \cite{Stembridge:B}. Elements of Coxeter groups give rise to heaps as follows.

\begin{definition}%3.1.5
Let {\it X} be a graph, let {\it S} be the set of vertices of {\it X} and let $\mathcal{C}$ be the relation on {\it S} defined by $s_i  \ \mathcal{C} \ s_j$ if and only if $s_i = s_j$ or $s_i$ and $s_j$ are adjacent vertices. Let $ w = s_{i_1}s_{i_2} \cdots s_{i_l}$ be an arbitrary word in the generators $S$ of a Coxeter group  $W=W(X)$. The word $w$ gives a labelled heap $(E, \leq_{E}, \varepsilon)$ where $ E = \{1, 2, \ldots, l \}, \ \varepsilon(j)= s_{i_j}$, and the relation $\leq_\mathcal{C}$ of condition (ii) of Definition 3.1.1 can be defined by
\begin{displaymath}
a \leq_\mathcal{C} b \Leftrightarrow a \leq b    \quad \textnormal{and} \quad  \varepsilon(a) \ \mathcal{C} \   \varepsilon(b),
\end{displaymath}
where $\leq$ is the usual ordering on integers. The partial order $\leq_E$ is the reflexive and transitive closure of $\leq_\mathcal{C}$, and the heap of the word {\it w} is by definition the heap $E_l$ corresponding to the given labelled heap.
\end{definition}

The following example should give some insight into the previous definitions.
\begin{example}%3.1.6
Let {\it X} be the Coxeter graph of type $\widetilde{C}_7$ as shown in Figure 3 below. Following the previous definition, $S = \{s_1, s_2, s_3, s_4, s_5, s_6, s_7, s_8\}$. Let $w= s_1s_3s_5s_2s_4s_6s_1s_3s_5s_7 = s_{i_1} s_{i_2} \cdots s_{i_9}s_{i_{10}}$. The word {\it w} gives a labelled heap $E_{10}$ where $ E = \{1, 2, \ldots, 10 \}$, and $\varepsilon(j)= s_{i_j}$. Note that $s_{i_1}$ is the label for the vertex $s_1 \in X$ and $s_{i_4}$ is the label for the vertex $s_2 \in X$. The vertices $s_1$ and $s_2$ are adjacent vertices in $X$ so $s_{i_1} \leq_\mathcal{C} s_{i_4}$. The heap is shown in Figure 3(a).
\end{example}

\begin{center}

\includegraphics{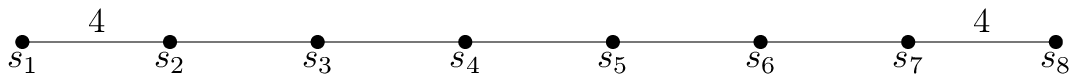}

FIG. 3. Coxeter graph of type $\widetilde{C}_7$.
\end{center}
\vspace{0.5 in}

\begin{center}

\includegraphics{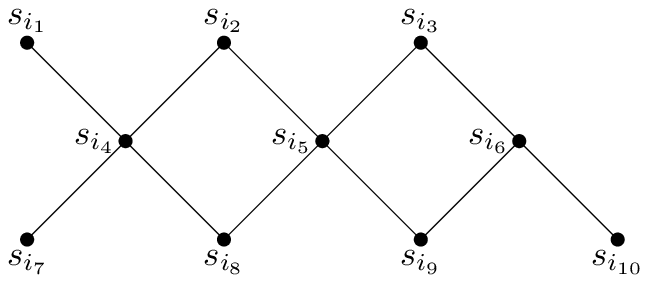}

FIG 3(a). The heap of $w = s_1s_3s_5s_2s_4s_6s_1s_3s_5s_7$ in type $\widetilde{C}_7$.
\end{center}

\begin{definition}%3.1.7
Let $(E, \leq, \varepsilon)$ be a labelled heap with pieces in {\it P} and let {\it F} be a subset of {\it E}. Let $\varepsilon'$ be the restriction of $\varepsilon$ to $F$. Let $\mathcal{R}$ be the relation defined on $F$ by $a \ \mathcal{R} \ b$ if and only if $a \leq b$ and  $\varepsilon{(a)} \ \mathcal{C} \ \varepsilon{(b)}.$ Let $\leq'$ be the transitive closure of $\mathcal{R}$. Then $(F, \leq', \varepsilon')$ is a labelled heap with pieces in $P$. The heap $[F, \leq', \varepsilon']$ is called a $subheap$ of $[E, \leq, \varepsilon]$. If {\it F} is convex as a subset of {\it E} (if $x_i \leq x_j \leq x_k \ \textrm{with} \ x_i, x_k \in F, \ \textrm{then} \ x_j \in F)$, we call {\it F} a {\it convex subheap} of {\it E}. A {\it trivial heap} is a heap $[E, \leq, \varepsilon]$ for which the order relation $\leq$ is trivial, meaning that no element of $E$ covers any other element.
\end{definition}

\begin{definition}%3.1.8
Let $E = [E, \leq_E, \varepsilon]$ and $F = [F, \leq_F, \varepsilon']$ be two heaps in $H(P, \mathcal{C})$. We define the heap $G = [G, \leq_G, \varepsilon''] = E \odot F$ of $H(P, \mathcal{C})$ (which we call superposition of $E$ over $F$ as follows.\\
(1) The underlying set $G$ is the disjoint union of $E$ and $F$.\\
(2) The labelling map $\varepsilon''$ is the unique map $\varepsilon'' : G \longrightarrow P$ whose restriction to $E$ (respectively, $F$) is $\varepsilon$ (respectively, $\varepsilon'$).\\
(3) The order relation $\leq_G$ is the transitive closure of the relation $\mathcal{R}$ on $G$, where $a \ \mathcal{R} \ b$ if and only if one of the following three conditions holds:\\
(i) $a, b \in E \ \textnormal{and} \ a \leq_E b;$\\
(ii) $a, b \in F \ \textnormal{and} \ a \leq_F b;$\\
(iii) $a \in E, b \in F \ \textnormal{and} \ \varepsilon(a) \ \mathcal{C} \ \varepsilon'(b).$
\end{definition}

\section{The Boundary Map $\partial$}
The following analogue to a certain boundary map $\partial$ in algebraic topology comes from \cite{Green:B}. 

\begin{definition}%3.2.1
Let $k$ be a field. Let $V_0$ be the set of elements of $[E, \leq, \varepsilon]$, i.e., the set of elements of (a representative of) the underlying poset, $E$.  We call the elements of  $V_0$ {\it vertices} and denote their {\it k}-span by $C_0$.

Let $V_1$ be the set of all pairs $(x,y) \in E \times E$ with $x<y$ and $\varepsilon(x) = \varepsilon(y)$ such that there is no element $z$ for which we have both $\varepsilon(x) = \varepsilon(z) = \varepsilon(y)$ and $ x < z < y$. We call the elements of $V_1$ {\it edges} and denote their {\it k}-span by $C_1$.

For all other integers $i \in \mathbb{Z} \backslash \{0, 1\}$, we define $C_i = 0$.

The {\it k}-linear map $\partial = \partial_E : C_1 \rightarrow C_0$ is defined by its effect on the edges as follows:
\begin{displaymath}
\partial : (x,y) \mapsto  \sum_{\substack{x<w<y \\ \varepsilon(w) \ \mathcal{C} \ \varepsilon(x)}} w.
\end{displaymath}
\end{definition}

\begin{example}%3.2.2
Consider the heap {\it E} arising from Example 3.1.6. The edges are $e_1 = (1, 7)$, $e_2 = (2, 8)$, and $e_3 = (3, 9)$. However, we will abuse notation and write $s_{i_k} = k$ so that the sums we create using the map $\partial$ make sense. Thus, $e_1 = (s_{i_1}, s_{i_7})$, $e_2 = (s_{i_2}, s_{i_8})$, and $e_3 = (s_{i_3}, s_{i_9})$. We now have $\partial(e_1) = s_{i_4}$, $\partial(e_2) = s_{i_4} + s_{i_5}$, and $\partial(e_3) = s_{i_5} + s_{i_6}$.
\end{example}

\begin{definition}%3.2.3
Let $E_n = E$ be a heap and let {\it k} be a field. If $ v \in E$, we let $E(v)$ be the subheap of {\it E} obtained by defining $E(v) = E \backslash \{v\}$. We say {\it E} is {\it acyclic} if ker $\partial_E=0$. We say {\it E} is {\it strongly acyclic} if {\it E} is acyclic and $E(v)$ is acyclic for all $ v \in E$. We say {\it v} is a {\it boundary vertex} of {\it E} if $v \in \textrm{Im}(\partial_E)$. We say {\it v} is an {\it effective boundary vertex} if $\partial_{E}(e_0) = v$ for an edge $e_0$ in $E_n$. Let $R$ be the relation defined on the vertices of $E_n$ by $v_1 \ R \ v_2$ if and only if  $\partial_{E}(e_0) = v_1 + v_2$ for some edge $e_0$. Linear equivalence is the reflexive, transitive closure of $R$. Since $R$ is also symmetric, linear equivalence is an equivalence relation on the vertices of $E$. We will sometimes write $R_E$ instead of $R$ to show the relation holds on a particular set $E$.
\end{definition}

\begin{remark}%3.2.4
The definitions for linear equivalence and effective boundary vertices are inspired by the notation and ideas of the Riemann--Roch theorem for graphs \cite{Baker;S.Norine:A}.
\end{remark}
\begin{example}%3.2.5
In Example 3.1.6, since $\partial(e_1) = s_{i_4}$, we see that $s_{i_4}$ is an effective boundary vertex. Furthermore, since $\partial(e_2) = s_{i_4} + s_{i_5}$ and $\partial(e_3) = s_{i_5} + s_{i_6}$, $s_{i_4}$ and $s_{i_5}$ are linearly equivalent and $s_{i_5}$ and $s_{i_6}$ are linearly equivalent. By the transitivity of linear equivalence, $s_{i_6}$ and $s_{i_4}$ are linearly equivalent. The equivalence classes are $\{s_{i_4}, s_{i_5}, s_{i_6}\}$, and the singletons $\{s_{i_1}\}$, $\{s_{i_2}\}$, $\{s_{i_3}\}$, $\{s_{i_7}\}$, $\{s_{i_8}\}$, $\{s_{i_9}\}$, $\{s_{i_{10}}\}$. It is clear that dim $\textrm{Im}(\partial_E) = 3.$ There are 3 edges, therefore ker $\partial_E=0$ and {\it E} is acyclic. If we consider the subheap $E(s_{i_5})$,  we find that $\partial(e_1) = s_{i_4}$ and $\partial(e_2) = s_{i_4}$, thus ker $\partial_{E(s_{i_5})} \neq 0$; hence the heap {\it E} is not strongly acyclic. Note that $\partial_E(e_2 - e_1) = (s_{i_4} + s_{i_5}) - (s_{i_4}) = s_{i_5}$, thus $s_{i_5}$ is a boundary vertex; similarly, $s_{i_6}$ is a boundary vertex. 
\end{example}

\section{Some Properties of Heaps and Full Commutativity}
We recall property P1 from \cite{Green:B}. 

\begin{definition}(\emph{Property P1})%3.3.1
\ Let $E= [E, \leq, \varepsilon]$ be a heap. We write $E(a) \prec^{+} E$ (respectively, $E(a) \prec^{-} E$) if $a$ is a maximal (respectively, minimal) vertex of $E$ and there exists a maximal (respectively, minimal) vertex $b$ of $E(a)$ with $\varepsilon(b) \neq \varepsilon(a)$ such that $b$ is not maximal (respectively, minimal) in $E$. We write  $E(a) \prec E$ if either  $E(a) \prec^{+} E$ or  $E(a) \prec^{-} E$.
If there is a (possibly trivial) sequence $E_1 \prec E_2 \prec \cdots \prec E$ of heaps in $H(P, \mathcal{C})$ where $E_1$ is a trivial heap, we say that the heap {\it E} is {\it dismantlable} or that {\it E} has property P1.
\end{definition}

\begin{example}%3.3.2
The heap arising from Example 3.1.6 is dismantlable and the following is a suitable chain. Note that $E(v_1, \ldots, v_k)$ is the subheap corresponding to $E \backslash \{v_1, \ldots, v_k\}$.
\begin{align*}
E(s_{i_4}, s_{i_5}, s_{i_6}, s_{i_1}, s_{i_2}, s_{i_3}) &\prec E(s_{i_5}, s_{i_6}, s_{i_1}, s_{i_2}, s_{i_3})\\
&\prec E(s_{i_6}, s_{i_1}, s_{i_2}, s_{i_3})\\
 &\prec E(s_{i_1}, s_{i_2}, s_{i_3}) \\
 &\prec E(s_{i_2}, s_{i_3})\\
 &\prec E(s_{i_3}) \prec E
\end{align*}
where $E_1 =  E(s_{i_4}, s_{i_5}, s_{i_6}, s_{i_1}, s_{i_2}, s_{i_3})$ is a trivial heap.
\end{example}

Let $(W, S)$ be a Coxeter system with Coxeter graph $X$. Denote by $S^*$ the free monoid on $S$. We call the elements of {\it S letters} and those of $S^*$ {\it words}. Let $\phi_X : S^* \longrightarrow  W$ be the surjective morphism of monoid structures satisfying $\phi_X(i) = i$ for all $i \in S$. A word \emph{i} $\in S^*$ is said to {\it represent} its image $w = \phi_X(\emph{i}) \in W$; furthermore, if the length of \emph{i} is minimal among the lengths of all the words that represent {\it w}, then \emph{i} is a reduced expression for {\it w}. The {\it commutation monoid} Co$(X, S)$ is the quotient of the free monoid $S^*$ by the congruence $\equiv$ generated by the commutation relations:
\begin{displaymath}
st \equiv ts \ \textnormal{for all}\ s, t \in S \ \textnormal{with} \ \phi_X(s)\phi_X(t) = \phi_X(t)\phi_X(s).
\end{displaymath}
Note that, as a monoid, {\it W} is a quotient of Co$(X, S)$. We define $\boldsymbol{[i]}$ to be the image of $\boldsymbol{i} \in S^*$ in Co$(X, S)$. We will sometimes refer to the elements of Co$(X, S)$ as $traces$. 

The product $\odot$ of heaps from Definition 3.1.8 is associative and $H(S, \mathcal{C})$ is a $monoid$, called the $heap \ monoid$, whose identity element is the empty heap.

We define the map $\phi_V : S^* \longrightarrow  H(S, \mathcal{C})$ by the relation
\begin{displaymath}
\textnormal{for} \ w = \alpha_1\alpha_2 \cdots \alpha_n \in S^*, \ \phi_V(w) = \alpha_1 \odot \alpha_2 \odot \cdots \odot \alpha_n \in H(S, \mathcal{C})
\end{displaymath}

\begin{prop}%3.3.3
Let $H(S, \mathcal{C})$ be a heap monoid with pieces in $S$ and concurrency relation $\mathcal{C}$.
Let $C$ be the complementary relation of $\mathcal{C}$. The morphism of monoids $\phi_V : S^*  \longrightarrow  H(S, \mathcal{C})$ defined above induces an isomorphism 
$\overline{\phi}$ 
between the monoid $H(S, \mathcal{C})$ and the commutation monoid Co($X, S)$.
\end{prop}
\emph{Proof.} This is a restatement of \cite[Proposition 3.1]{Viennot:A}. \QED

\begin{lemma}%3.3.4
Let Co($X, S$) be the commutation monoid defined above. Then the map $\psi: \textnormal{Co}(X, S)  \longrightarrow  W$ is a surjective morphism of monoid structures satisfying $\psi(\boldsymbol{[s]}) = \phi_X(\boldsymbol{s})$. 
\end{lemma}
\emph{Proof.} The relations in Co($X, S$) are all relations in $W$, hence the map is well defined. \QED

\begin{theorem}%3.3.5
Maintain the above notation. The following diagram commutes, where all the maps are morphisms of monoids.
\end{theorem}
\begin{displaymath} 
\xymatrix{ 
S^* \ar[r]^{\phi_V} \ar[dr]_{\phi_X} & \textnormal{Co}(X, S) \ar[d]^{\psi} \ar@{<->}[r]^{\overline{\phi}} & H(S, \mathcal{C}) \\ 
                                                             & W  } 
\end{displaymath}
\emph{Proof.} 
The diagram collects the results from Proposition 3.3.3, Lemma 3.3.4, and the definitions preceding Proposition 3.3.3. \QED

\begin{definition}%3.3.6
The equivalence class $\boldsymbol{[s]} = \phi_V(\boldsymbol{s})$ of a given reduced word $\boldsymbol{s} \in S^*$ (relative to $\equiv$) consists of the words obtainable from $\boldsymbol{s}$ by transposing adjacent commuting pairs of generators from $S$. We call $\boldsymbol{[s]}$ the {\it commutativity class} of $\boldsymbol{s}$. If we consider the set $R(w) \subset S^*$ of all reduced expressions for some $w \in W$, $R(w)$ is a union of commutativity classes, namely the $\phi_V(R(w))$. If $\phi_V(R(w))$ is a singleton, we say that $w$ is fully commutative. Note that this is a characterization of full commutativity. If $w \in W_c$, and $\boldsymbol{s}$ is a reduced expression for $w$, then $\phi_V(\boldsymbol{s})$ is independent of the choice of $\boldsymbol{s}$. We call $\overline{\phi}(\phi_V(\boldsymbol{s}))$ the {\it heap of w}.
\end{definition}

\begin{definition}%3.3.7
Let $E$ be a poset with partial order $\leq$. By a $subposet$ of $E$, we mean a subset $X$ of $E$ equipped with the induced partial order; that is, the partial ordering of $X$ such that for $x, y \in X$ we have $x \leq y$ in $X$ if and only if $x \leq y$ in $E$. We define the {\it closed interval} $[x, y]_E$ to be $\{z \in E : x \leq z \leq y\}$. If $X$ is totally ordered, then $X$ is called a $chain$. A subposet $X$ of $E$ is $convex$ if $y \in X$ whenever $x < y < z$ in $E$ and $x, z \in X$. (In particular, an interval is convex.)
\end{definition}

Suppose that {\it E} is the heap of some $w \in W_c$. For each $s \in S$, the members of {\it E} with label {\it s} form a chain. We denote by $s^{(k)}$ the $k$-th least member of this chain with respect to {\it E}.

Some of the main results from \cite{Stembridge:B} are summarized in the following theorem.

\begin{theorem}  % 3.3.8
(i) The heap {\it E} of a word $\emph{s} \in S^*$ is the heap of some fully commutative element $w \in W$ if and only if the following two conditions hold:\\
\indent (a) there is no convex chain $i_1 < \cdots < i_m$ in $E$ such that $s_{i_1} = s_{i_3} = \cdots = s$ \indent and $s_{i_2} = s_{i_4} = \cdots = t$, where $3 \leq m = m(s, t) < \infty$, and \\
\indent (b) there is no covering relation $i < j$ in $E$ such that $s_i = s_j$.\\
(ii) Subwords of fully commutative elements are fully commutative.\\
(iii) Let $s \in S$, and let $w \in W$ be fully commutative with heap $E$. If $sw$ is not fully commutative, then $sw$ is reduced and there is a unique $t \in S$ such that $m(s, t) \geq 3$ and $t^{(1)} < s^{(1)}$ in
{\it E}. Moreover, $m(s, t) < \infty$ and
\begin{displaymath}
 t^{(1)} < s^{(1)} < t^{(2)} < \cdots < s^{(k)} \quad (\textrm{if} \ m(s, t) \ = 2k+1)
\end{displaymath}
\begin{displaymath}
t^{(1)} < s^{(1)} < t^{(2)} < \cdots  < t^{(k)} \quad (\textrm{if} \ m(s, t) \ = 2k),
\end{displaymath}
is a convex chain in {\it E}.
\end{theorem}
\emph{Proof.} For (i), see \cite[Proposition 2.3]{Stembridge:B}. For (ii), see \cite[Prop 1.1]{Stembridge:B}.
For (iii), see \cite[Lemma 3.1]{Stembridge:B}. \QED

\begin{lemma}%3.3.9
Let $E$ be a finite poset. Then $X \subset E$ is a convex subposet of $E$ if and only if there is a sequence $E_0, E_1, \ldots, E_k$ of subposets of $E$ such that
\begin{displaymath}
X = E_k \subset E_{k-1} \subset \cdots \subset E_1 \subset E_0 = E
\end{displaymath}
where $E_i \backslash E_{i + 1} = \{\alpha\}$ with $\alpha$ maximal or minimal in $E_i$.
\end{lemma}
\emph{Proof.} Assume that $X$ is a convex subposet of $E$. The proof is by induction on $k$, where $k = |E \backslash X|$. If $k = 0$, then $X = E$ and there is nothing to prove. Suppose that $k > 0$. Since $X$ is convex by assumption, $X$ cannot contain all the maximal and minimal elements of $E$. Thus we can pick a maximal (or minimal) element $\alpha \in E$ such that $\alpha \notin X$. Since $X$ is a convex subset of $E \backslash \{\alpha\}$ and $|(E \backslash \{\alpha\}) \backslash X| = k-1$, we can apply the induction hypothesis. The conclusion follows.

Conversely, assume that there is a sequence of subposets of $E$ such that
\begin{displaymath}
X = E_k \subset E_{k-1} \subset \cdots \subset E_1 \subset E_0 = E
\end{displaymath}
where $E_i \backslash E_{i + 1} = \{\alpha\}$ with $\alpha$ maximal or minimal in $E_i$. Note that ``is a convex subset of " is a transitive relation of subposets of a poset $E$. If we remove a maximal or minimal vertex $\alpha \in E$, we have a convex subset $E \backslash \{\alpha\}$. Combining these two facts gives us the result. \QED

\begin{lemma}%3.3.10
Let $X$ be a convex subheap of a heap $E$. If $e_0$ is an edge in $X$, then $e_0$ is an edge in $E$.
\end{lemma}
 \emph{Proof.} By way of contradiction, let $e_0 = (x, y)$ be an edge in $X$ that is not an edge in $E$. Then there exists $z \in E$ such that $x < z < y$ and $\varepsilon(x) = \varepsilon(z) = \varepsilon(y)$, and $z$ must be in $X$ since $X$ is a convex subset of $E$. Therefore $e_0$ is an edge in $E$. \QED
 
 \begin{lemma}%3.3.11
 Let $E$ be a finite poset with convex subposet $X$. Then $[x, y]_X = [x, y]_E$.
 \end{lemma}
 \emph{Proof.} 
 Let $z \in [x,y]_X$. Then $z \in X$ and $x \leq z \leq y$. Thus $z \in [x, y]_E$. Conversely, let $z \in [x, y]_E$ with $x \leq z \leq y$. Since $x, y \in X$, by convexity we have $z \in X$. Thus $z \in [x, y]_X$. \QED
 
 \begin{lemma}%3.3.12
Let $E$ be a heap and let $X$ be a convex subheap of $E$. Let $e = (x, y)$ be an arbitrary edge in $X$. Identifying $X$ as a subset of $E$, we have $\partial_E(e) = \partial_X(e)$.
\end{lemma}
\emph{Proof.}
By definition of $\partial_X$, we have 
\begin{align*}
 \partial_X(e) &=  \sum_{\substack{x<z<y \\ \varepsilon(z) \ \mathcal{C} \ \varepsilon(x) \\ z \in X}}z \ =  \sum_{\substack{x<z<y \\ \varepsilon(z) \ \mathcal{C} \ \varepsilon(x) \\ z \in E}}z\\
   &= \partial_E(e),
 \end{align*}
  where the equality of sums is by Lemma 3.3.11. \QED

\begin{lemma}%3.3.13
Let $E$ be a heap and let $X$ be a convex subheap of $E$.\\
(i) If $\alpha, \beta$ are linearly equivalent vertices in $X$, then $\alpha, \beta$ are linearly equivalent in $E$.\\
(ii) If $\alpha$ is an effective boundary vertex in $X$, then $\alpha$ is an effective boundary vertex in $E$.
\end{lemma}
\emph{Proof.}
For (i), it is enough to assume that $\alpha \ R_X \ \beta$ in $X$. Then there is an edge of $X$, say $e_0$, with $\partial_X(e_0) = \alpha + \beta$. By Lemma 3.3.10, $e_0$ is an edge of $E$. By Lemma 3.3.12, $\partial_E(e_0) = \alpha + \beta$. Therefore $\alpha \ R_E \ \beta$. \\
For (ii), let $\alpha$ be an effective boundary vertex in $X$. Then there is an edge $e_0$ in $X$ with $\partial_X(e_0) = \alpha$. By Lemma 3.3.10, $e_0$ is an edge of $E$. By Lemma 3.3.12,  $\partial_E(e_0) = \alpha$. Thus $\alpha$ is an effective boundary vertex of $E$. \QED

\begin{lemma}%3.3.14
Let $W$ be a star reducible Coxeter group and let $w \in W_c$. Let $E$ be the heap of $w$ and let $X$ be a convex subheap. Let $B$ and $B'$ be the sets of boundary vertices in $E$ and $X$, respectively. Then $B' \subset B.$
\end{lemma}
\emph{Proof.}
Assume for a contradiction that $\alpha \in B'$ and $\alpha \notin B$. By definition, there exist edges $e_0, e_1, \ldots, e_r$ in $X$ such that $$\partial_X \left( \sum_{i=0}^r\lambda_i e_i \right) = \alpha.$$ By Lemma 3.3.12, $\partial_E(e_i) = \partial_X(e_i)$ for each $i = 0, 1, \ldots, r$. Thus $$\partial_X \left( \sum_{i=0}^r\lambda_i e_i \right) = \partial_E \left( \sum_{i=0}^r\lambda_i e_i \right) = \alpha. \qquad \blacksquare$$ 

\begin{lemma}%3.3.15
Let $E$ be an arbitrary heap. If $\alpha, \beta$ are linearly equivalent in $E$ and $\alpha$ is a boundary vertex, then $\beta$ is a boundary vertex in $E$.
\end{lemma}
\emph{Proof.}
We may assume that $\alpha \ R \ \beta$. By hypothesis, if $e_0, e_1, \ldots, e_r$ are the edges of $E$, then $\partial_E(e_k) = \alpha + \beta$ for some $k$, and $$\partial_E \left( \sum_{i=0}^r\lambda_i e_i \right) = \alpha$$ for some scalars $\lambda_1, \lambda_2, \ldots, \lambda_r$. Hence $$ \beta  = \partial(e_k) - \alpha = \partial(e_k) - \partial \left( \sum_{i=0}^r\lambda_i e_i \right) = \partial\left(e_k - \sum_{i=1}^r\lambda_i e_i \right),$$ which shows that $\beta$ is a boundary vertex.
 \QED
 
 \begin{definition} \cite[Property P2]{Green:B} \ %3.3.16
 We say a heap $E$ has property P2 if it contains no convex chains of the form $x < y < z$ or $x < z$ with $\varepsilon(x) = \varepsilon(z)$ in either case.
\end{definition}

\begin{example}%3.3.17
The heap arising from Example 3.1.6 does not have property P2. Although there are no convex chains of the form $s_{i_k} <  s_{i_j}$ with $\varepsilon(s_{i_k}) = \varepsilon(s_{i_j})$, the chain $s_{i_1} <  s_{i_4} < s_{i_7}$ violates the other requirement.
\end{example}

The following theorem relates star reducible Coxeter groups, properties P1 and P2, and the definitions for strongly acyclic and acyclic heaps.
\begin{theorem}%3.3.18
Let $W = W(X)$ be a star reducible Coxeter group with Coxeter graph $X$ and suppose that $E$ is the heap of (a reduced word for) some $w \in W_c$. Then\\
(i) {\it E} has property P1 if and only if {\it E} is acyclic;\\
(ii) if {\it E} is strongly acyclic, then $E$ has property P2;\\
(iii) if $X$ is simply laced, then the converse to (ii) holds;\\
(iv) any heap for $w \in W_c$ has property P1.
\end{theorem}
\emph{Proof.} Part (i) is \cite[Theorem 2.4.4]{Green:B}. Part (ii) is \cite[Proposition 2.2.7]{Green:B}. Part (iii) is \cite[Theorem 2.4.2 (i)]{Green:B}. Part (iv) follows directly from the definition of star reducibility. \QED

\begin{lemma} %3.3.19
Let $X$ be a bipartite, simply laced Coxeter graph with $X = X_1 \dot\cup X_2$, where $X_1$ and $X_2$ each consists of mutually nonadjacent vertices. Let $W = W(X)$ be the corresponding Coxeter group. Let $E$ be the heap of some word on the Coxeter generators of $W$. \\
(i) If $e= (x, y)$ and $\partial_E(e) = \sum_{j=0}^r\lambda_j \alpha_j$, then $x, y \in X_i$ for some $i \in \{1, 2\}$ and all $\alpha_j$ with $\lambda_j \neq 0$ lie in $X_{3-i}$.\\
(ii) If $E$ is acyclic and $\alpha$ is a boundary vertex of $E$, then $E \backslash \{\alpha\}$ is not acyclic.\\
(iii) If $E$ is the heap of a fully commutative element of $W$, then $E$ has no boundary vertices.
\end{lemma}
\emph{Proof.}
For (i), note that
\begin{displaymath}
\partial((x,y)) =  \sum_{\substack{x< z <y \\ \varepsilon(x) \ \mathcal{C} \ \varepsilon(z)}} z.
\end{displaymath}
Since $(x, y)$ is an edge, $\varepsilon(z)$ is adjacent to $\varepsilon(x)$. Therefore, if $\varepsilon(x) = \varepsilon(y) \in X_i$, then $\varepsilon(z) \in X_{3 - i}$. The assertion follows.

For (ii), assume that $\alpha$ is a boundary vertex of $E$. By definition, there exist edges $e_0, e_1, \ldots, e_r$ in $E$ such that $$\partial_E \left( \sum_{i=0}^r\lambda_i e_i \right) = \alpha.$$ By (i), $\alpha$ is not an endpoint of any edge $e_i$, so deleting $\alpha$ gives us
$$\partial_{E \backslash \{\alpha\}} \left( \sum_{i=0}^r\lambda_i e_i \right) = 0.$$ But then ker $\partial_{E \backslash \{\alpha\}} \neq 0$ and hence $E \backslash \{\alpha\}$ is not acyclic.

For (iii), assume by way of contradiction that $\alpha$ is a boundary vertex in $E$. By definition, there exist edges $e_0, e_1, \ldots, e_r$ in $E$ such that $$\partial_E \left( \sum_{i=0}^r\lambda_i e_i \right) = \alpha.$$ By (i), $\alpha$ is not an endpoint of any edge $e_i$. Since $X$ is simply laced, $E$ has Property P2 and by Theorem 3.3.18 (iii), $E$ is strongly acyclic. By (ii), if we delete the vertex $\alpha$, we have
that $E \backslash \{\alpha\}$ is not acyclic and hence $E$ is not strongly acyclic, a contradiction.
 \QED

 \section{Main Theorem}

\begin{theorem}%3.4.1
Let $X$ be a bipartite Coxeter graph with star reducible Coxeter group $W = W(X)$. In the heap $E$ of a fully commutative element $w \in W$, every boundary vertex is linearly equivalent to an effective boundary vertex.
\end{theorem}
\emph{Proof.} If $X$ is simply laced, then the claim is true for $W(X)$ by Lemma 3.3.19 (iii). By the classification of star reducible Coxeter groups in \cite{Green:P}, we can assume $X$ is a straight line graph. Let $x \in W_c$, and let $E$ be the heap of $x$ as in Definition 3.3.6.

The proof is by induction on $n$, the number of vertices in the heap $E$. If $x$ is a product of commuting generators from {\it S}, then the heap of $x$ has no edges and hence no boundary vertices. Therefore the claim is true vacuously and the cases $n = 0$ and $n=1$ are covered.

Since $W$ is a star reducible Coxeter group, we may assume that $x$ either has a reduced expression beginning with $st$ or ending with $ts$, where $I = \{s, t\}$ is a pair of noncommuting generators in $S$. We deal here with the case where $x$ has a reduced expression beginning with $st$. If $x$ has a reduced expression ending in $ts$, a symmetrical argument gives the same conclusion.

Taking star operations with respect to $I$, we can star reduce $x = stv$ (reduced) to $w=tv$. Since $w$ is a subword of the fully commutative element $x$, by Theorem 3.3.8 (ii), $w$ is fully commutative. Since the heap $E_{n-1}$ of $w$ is a convex subheap of $E_n$, by Lemma 3.3.14, the boundary vertices $B'$ of $E_{n-1}$ are a subset of the boundary vertices $B$ of $E_n$. If $B' = B$, by Lemma 3.3.13 the claim is true.

Suppose there is a boundary vertex $\alpha \in B$ with $\alpha \notin B'$. Then there exist edges $e_0, e_1, \ldots e_r$ in $E_n$ such that 
\begin{displaymath}
\partial_E \left( \sum_{i=0}^r\lambda_i e_i \right) = \alpha. 
\end{displaymath}
where $e_0 = (s^{(0)}, s^{(1)})$ with $s^{(0)} = s$ in the definition of $x$. Because $t^{(1)}$ is minimal in $E \backslash \{s^{(0)}\}$, $t^{(1)}$ does not appear in $\partial_E(e_i)$ for any $i \neq 0.$ Thus $t^{(1)} = \alpha$, $\lambda_0 = 1$, and $t^{(1)}$ is unique. We want to show that $t^{(1)}$ is linearly equivalent to a boundary vertex $\beta \in B'$. 

Since $t^{(1)}$ is a boundary vertex in $E_n$, $E_n(t^{(1)})$ has a cycle by Lemma 3.3.19 (ii). By Theorem 3.3.18 (i), $E_n(t^{(1)})$ does not have Property P1. Therefore $E_n(t^{(1)})$ is not the heap of a fully commutative element by Theorem 3.3.18 (iv). By Theorem 3.3.8 (i), either $E_n(t^{(1)})$ has a convex chain $s^{(0)} < u^{(1)} < s^{(1)} < \cdots < u^{(k)}$ (or $s^{(k)}$, depending on $m(s, u)$) or a convex chain $s^{(0)} < s^{(1)}$. If the second statement is true, then $s^{(0)} < t^{(1)} < s^{(1)}$ in $E_n$, $t^{(1)}$ is an effective vertex in $E_n$ and we are done. 

If the first statement is true, then we must have $\partial_E(s^{(0)}, s^{(1)}) = u^{(1)} + t^{(1)}$ and $u^{(1)}$ is a boundary vertex in $E_n$, by Lemma 3.3.15. As $t^{(1)}$ is the only boundary vertex in $B$ that is not in $B'$, $u^{(1)}$ must be a boundary vertex in $B'$. By induction, $u^{(1)}$ is linearly equivalent in $E_{n-1}$ to an effective boundary vertex in $E_{n-1}$ and thus in $E_n$ by Lemma 3.3.13. As $t^{(1)}$ is linearly equivalent to $u^{(1)}$ in $E_n$, by the transitivity of linear equivalence, $t^{(1)}$ is linearly equivalent in $E_n$ to an effective boundary vertex in $E_n$, and the claim is true. \QED 

The following is an example to illustrate that the theorem holds.

\begin{example}
Let $X$ be the Coxeter graph of type $\widetilde{C}_7$, as shown in Figure 5 below. Since $X$ has no odd cycles, we see that $X$ is bipartite. Also note that $W(X)$ is star reducible by the classification in \cite{Green:P}. Let $w= s_1s_3s_5s_2s_4s_6s_1s_3s_5s_7 = s_{i_1} s_{i_2} \cdots s_{i_9}s_{i_{10}}$ and let $E$ be the heap of $w$. In Example 3.2.5, we showed that the boundary vertices of $E$ are $s_{i_4}, s_{i_5}, s_{i_6}$. Because  $\partial(e_1) = s_{i_4}$, we see that $s_{i_4}$ is an effective boundary vertex. Furthermore, since $\partial(e_2) = s_{i_4} + s_{i_5}$ and $\partial(e_3) = s_{i_5} + s_{i_6}$, $s_{i_4}$ and $s_{i_5}$ are linearly equivalent and $s_{i_5}$ and $s_{i_6}$ are linearly equivalent. By the transitivity of linear equivalence, $s_{i_6}$ and $s_{i_4}$ are linearly equivalent. Thus $s_{i_5}$ and $s_{i_6}$ are both linearly equivalent to the effective boundary vertex $s_{i_4}$ and the theorem holds.
\end{example}

\begin{center}

\includegraphics{C7.mps.ps}

FIG. 5. Coxeter graph of type $\widetilde{C}_7$.
\end{center}
\vspace{0.5 in}

\begin{center}

\includegraphics{C7Heap.mps.ps}

FIG 5(a). The heap of $w = s_1s_3s_5s_2s_4s_6s_1s_3s_5s_7$ in type $\widetilde{C}_7$.
\end{center}

The following is a non-example.
\begin{example}
 Consider the bipartite Coxeter graph of Figure 6.
\begin{center}

\includegraphics{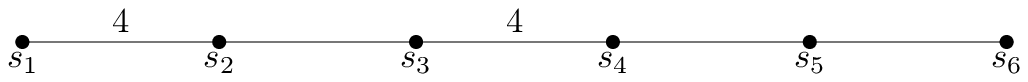}

FIG. 6. Coxeter graph on six vertices.
\end{center} 

Let $S = \{s_1, s_2, \ldots, s_6\}$. Let $w= s_1s_3s_5s_2s_4s_1s_3s_2s_4s_6s_1s_3s_5 = s_{i_1} s_{i_2} \cdots s_{i_{12}}s_{i_{13}}$. The word $w$ gives a labelled heap $E_{13}$ where $ E = \{1, 2, \ldots, 13 \}, \ \varepsilon(j)= s_{i_j}$. The heap is shown in Figure 7.
\end{example}

\begin{center}

\includegraphics{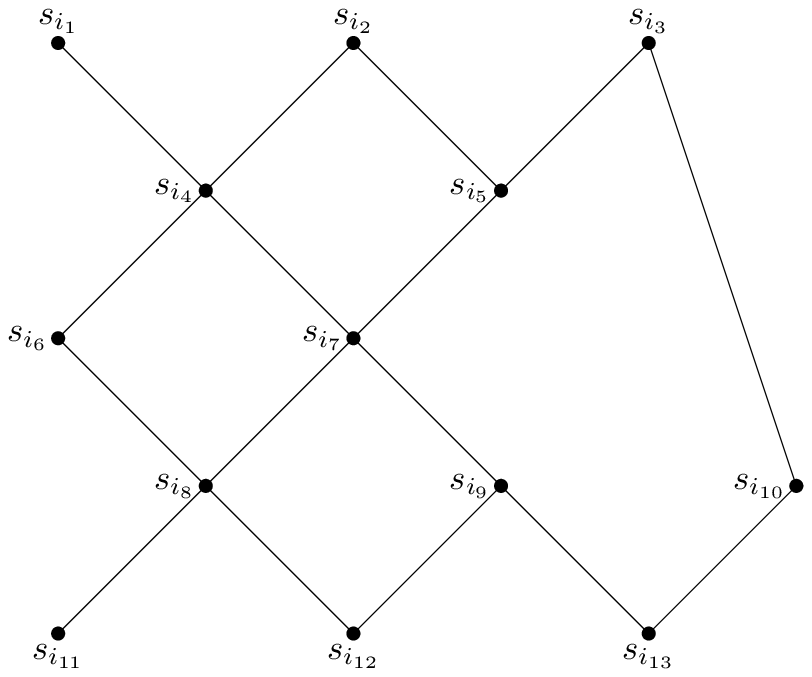}

FIG 7. The labelled heap of $w= s_1s_3s_5s_2s_4s_1s_3s_2s_4s_6s_1s_3s_5$.
\end{center}

Note that $w$  is fully commutative and the heap is dismantlable.
However $s_{i_{10}}$ is a boundary vertex that is not linearly equivalent to an effective boundary vertex. 
Theorem 3.4.1 does not apply in this case because the group defined by the Coxeter graph is not one of the star reducible Coxeter groups.

\chapter{Property W}
\section{Introduction} In this chapter we present an application of Theorem 3.4.1. 

 \begin{lemma} \label{lemma:L}
 Let $E$ be a heap and let $\alpha \in E$. Then $| \textrm{dim ker} \ \partial_{E} - \textrm{dim ker} \ \partial_{E(\alpha)}| \leq 1$
 \end{lemma}
 \emph{Proof.}
 This is \cite[Theorem 2.1.1]{Green:B}.
 
 The reader may wish to review the definitions of generalized Temperley--Lieb algebra and $\widetilde{t}$-basis given in \S{2.1}. 

 \begin{definition} \label{def:C}
 Let $W$ be a Coxeter group with Coxeter graph $X$ and let $w \in W_c$. Let $w = s_1s_2 \cdots s_r$ be a reduced expression for $w$. For each $s \in S$, let $b_s = v^{-1}\widetilde{t_1} + \widetilde{t}_s$, then define $b_w \in TL(X)$ by $b_w = b_{s_1}b_{s_2} \cdots b_{s_r}$. Note that the element $b_w$ is well-defined because any two reduced expressions for $w$ are commutation equivalent and $b_sb_t = b_tb_s$ if $st = ts.$ The set $\{b_w : w \in W_c\}$ is a basis for $TL(X)$, called the {\it monomial basis}; see Lemma~\ref{lemma:A} (ii) for more details.
 \end{definition}
 
 \begin{definition} \label{def:D}
 Let $W = W(X)$ be a star reducible Coxeter group. We define the function $h : S^* \longrightarrow \ZZ^{\geq 0}$ by $h(\boldsymbol{u}) = \textnormal{dim ker} \ \partial_{\overline{\phi}({\phi_V}(\boldsymbol{u}))}$. We sometimes write $h(\boldsymbol{[u]})$ instead of $h(\boldsymbol{u})$; this makes sense by Theorem 3.3.5.
 \end{definition}
 
  \begin{lemma} \label{lemma:S}
 Let $W(X)$ be a star reducible Coxeter group, let $\boldsymbol{s} = s_1s_2 \cdots s_r \in S^*$, and let $b(\boldsymbol{s}) \in TL(X)$ given by $b(\boldsymbol{s}) = b_{s_1}b_{s_2} \cdots b_{s_r}$. Express $b$ as a linear combination of the monomial basis, namely $$b(\boldsymbol{s}) = \sum_{w \in W_c} \lambda_wb_w.$$ Then each $\lambda_w$ is an integer multiple of $(v + v^{-1})^{h(\boldsymbol{s})}$.
\end{lemma}
\emph{Proof.}
This is \cite[Lemma 2.9]{Green:P}. \QED

Recall from \S{1} that $\mathcal{A^-}$ is the ring $\ZZ[v^{-1}]$.
 
 \begin{lemma}  
\label{lemma:A}
Let $W = W(X)$ be a star reducible Coxeter group with Coxeter graph $X$.\\
(i) If $x \in W$ is weakly complex, then $\widetilde{t}_x \in \mathcal{L}.$ \\
(ii) The $b$-basis and the $c$-basis of $TL(X)$ have the same $\ZZ$-span. In particular, the $b$-basis is an $\mathcal{A^-}$-basis for $\mathcal{L}$.
\end{lemma}
\emph{Proof.}
Part (i) is \cite[Lemma 3.8(i)]{Green:P}. Part (ii) is \cite[Lemma 2.10]{Green:P}. \QED

\begin{lemma} \label{lemma:M}
Let $\boldsymbol{s} = s_1s_2 \cdots s_r \in S^*$ with $\overline{\phi}(\phi_V(\boldsymbol{s}))$ acyclic and write
$$\widetilde{t}_{\boldsymbol{s}}  = \widetilde{t}_{s_1}\widetilde{t}_{s_2} \cdots \widetilde{t}_{s_r}.$$ Let $\widetilde{t}_{\boldsymbol{s'}}$ be the result of replacing $k$ distinct $\widetilde{t}_{s_i}$ with $(-v^{-1})$. Then $\widetilde{t}_{\boldsymbol{s'}} \in \mathcal{L}.$
\end{lemma}
\emph{Proof.}
Let $\boldsymbol{s} = s_1s_2 \cdots s_r$ be a reduced expression for $w \in W_c$ and write
\begin{align*}
\widetilde{t}_{\boldsymbol{s}} & = \widetilde{t}_{s_1}\widetilde{t}_{s_2} \cdots \widetilde{t}_{s_r}\\
& = (b_{s_1} - v^{-1})(b_{s_2} - v^{-1}) \cdots (b_{s_r} - v^{-1}).
\end{align*} Expanding the parentheses, we can express $\widetilde{t}_w$ as a linear combination of $2^r$ elements $(-v)^{-k}b(\boldsymbol{u})$, where \emph{u} is obtained from \emph{s} by deletion of $k$ generators. We have $$(-v)^{-k}b(\boldsymbol{u}) = \sum_{w \in W_c}(-v)^{-k}\lambda_wb_w.$$
By Lemma~\ref{lemma:S}, $\lambda_w \in v^{h(\UU)}\mathcal{A^-}$. By Lemma~\ref{lemma:L}, $h(\UU) \leq h(\ST) + k$. Since $\overline{\phi}(\phi_V(\ST))$ is acyclic by hypothesis, $h(\ST) = 0$. Thus $h(\UU) \leq k$ and $(-v^{-k})\lambda_w \in v^{h(\UU)-k}\mathcal{A^-} \subset \mathcal{A^-}$. By Lemma~\ref{lemma:A} (ii), $b_w \in \mathcal{L}$ so $(-v^{-k})\lambda_wb_w \in \mathcal{L}$ and $(-v^{-k})b(\UU) \in \mathcal{L}$.  

Now, $\widetilde{t}_{\boldsymbol{s'}}$ is a linear combination of $2^{r-k}$ terms of the form $(-v^{-k})b(\UU)$, and the conclusion follows. \QED

\begin{lemma} \label{lemma:Z}
Let $W$ be a Coxeter group and let $I = \{s, t\}$ be a pair of noncommuting generators in $S$. Suppose that whenever $x$ is weakly complex, $ux \in W_c$ and $u \in S$, we have $\widetilde{t}_x \in \mathcal{L}^u_L.$ Let $w = w_Iw^I$ be such that $w^I \in W_c$.\\
(i) If $sw_I < w_I$, then $\widetilde{t}_w \in \mathcal{L}^s_L$.\\
(ii) If $w_I = w_{st}$, the longest element in $W_I$, then $\widetilde{t}_{w_I}\widetilde{t}_{w^I} \in v^{-1}\mathcal{L}^s_L$.
\end{lemma}
\emph{Proof.} This is \cite[Lemma 4.11]{Green:K}. A right-handed analogue of this lemma also holds.\QED

\section{Acyclic Case}
\begin{lemma} \label{lemma:C}
Let $W(X)$ be a star reducible Coxeter group with a straight line Coxeter graph {\it X} and let $\boldsymbol{s}$ be a reduced expression for $w \in W_c$. If $x = sw \notin W_c$ and the heap $ E = \overline{\phi}(s\boldsymbol{[s]})$ is acyclic, then $\widetilde{t}_x \in v^{-1}\mathcal{L}$.
\end{lemma}
\emph{Proof.} Since $X$ is a straight line, each Coxeter generator fails to commute with at most two other generators. We will use this fact freely. The proof is by induction on $\ell(w)$. There are no weakly complex elements of length 0, thus the base case is vacuous.
By Theorem 3.3.8 (iii) it follows that $sw$ is reduced. From the proof of Theorem 3.3.8 (i), there exists a unique $t \in S$ such that $3 \leq m(s,t) < \infty$ and we have
\begin{displaymath}
s^{(1)} < t^{(1)} < s^{(2)} < t^{(2)} < \cdots < s^{(k + 1)} \quad (\textrm{if} \ m(s, t) \ = 2k+1)
\end{displaymath}
\begin{displaymath}
s^{(1)} < t^{(1)} < s^{(2)} < t^{(2)} < \cdots < s^{(k)} < t^{(k)} \quad (\textrm{if} \ m(s, t) \ = 2k),
\end{displaymath}
occurring as a convex chain \emph{c} in the heap $\overline{\phi}(s\boldsymbol{[s]})$, where $\boldsymbol{s}$ is a reduced expression for $w$. 

We first show that the five cases that follow are exhaustive. Since $E$ is acyclic, by Theorem 3.3.18 (i), $E$ has property P1. As $sw$ is weakly complex, $E$ is not trivial, so $E$ can be left or right star reduced. Thus there exist $\alpha, \beta \in E$ such that $E(\alpha) \prec^{+} E$ (respectively, $E(\alpha) \prec^{-} E$) if $\alpha$ is a maximal (respectively, minimal) vertex of $E$ and there exists a maximal (respectively, minimal) vertex $\beta$ of $E(\alpha)$ with $\varepsilon(\beta) \neq \varepsilon(\alpha)$ such that $\beta$ is not maximal (respectively, minimal) in $E$. This gives us two subcases for possible star reductions involving $\alpha$ and $\beta$. 

If $E(\alpha) \prec^{-} E$, then there exist elements $a, b \in S$ such that $a =\varepsilon(\alpha) \ \textnormal{and} \ b = \varepsilon(\beta)$. Either $\alpha \notin \boldsymbol{c}$ or $\alpha \in \boldsymbol{c}.$ If $\alpha \notin \boldsymbol{c}$, then $\beta \notin \boldsymbol{c}$ because $\beta$ is minimal in $E(\alpha)$. This is case (i) below. If $\alpha \in \boldsymbol{c}$, then $\alpha = s^{(1)}$. By the convexity of $\boldsymbol{c}$, $\beta = t^{(1)}$ and $\beta \in \boldsymbol{c}$; this is case (ii) below.

Suppose on the other hand that $E(\alpha) \prec^{+} E$. Either $\alpha, \beta \in \boldsymbol{c}, \alpha, \beta \notin \boldsymbol{c},$ or $\alpha \notin \boldsymbol{c}$ but $\beta \in \boldsymbol{c}.$ These are dealt with respectively in cases (iii), (iv), and (v) below.

\emph{Case (i)} In this case, both {\it a} and {\it b} commute with each of {\it s} and {\it t} and we have $w = abw'$ reduced. Since {\it x} is weakly complex, by Lemma~\ref{lemma:A} (i), $\widetilde{t}_x \in \mathcal{L}$ and $\ell(sw') > \ell(w')$. Also, since $x = sw = absw'$ is reduced, we can left star reduce {\it x} to the element $bsw'$. Now, $bsw'$ and $sw'$ are both weakly complex since $\alpha$ and $\beta$ are not involved in the chain \emph{c}. Thus the induction hypothesis is satisfied and we have $\widetilde{t}_{bsw'} \in v^{-1}\mathcal{L}$ and $\widetilde{t}_{sw'} \in v^{-1}\mathcal{L}$. Therefore $v\widetilde{t}_{sw'} \in \mathcal{L}$ and
\begin{displaymath}
\widetilde{t}_b(v\widetilde{t}_{sw'}) = v\widetilde{t}_{bsw'} \in \mathcal{L}^b
\end{displaymath}  
by Lemma~\ref{lemma:A} (i) and Proposition 2.1.9 (ii).
By Proposition 2.1.9 (iii), $\widetilde{t}_a(v\widetilde{t}_{bsw'}) \in \mathcal{L}^a$. Hence 
\begin{displaymath}
\widetilde{t}_a\widetilde{t}_{bsw'} = \widetilde{t}_{absw'} = \widetilde{t}_{x} \in v^{-1}\mathcal{L}^a \subset  v^{-1}\mathcal{L}.
\end{displaymath}  

\emph{Case (ii)} In this case, $a = s$, $b = t$, and we have $x = abw'$ reduced. We have $s\boldsymbol{[s]} = \boldsymbol{[cu]} \in \textnormal{Co}(X, S)$ with $\alpha, \beta \in \boldsymbol{c}$. Thus $sw = w_{st}w'$ reduced, where $w_{st}$ is the longest element in the parabolic subgroup generated by $I = \{s, t\} = \{a,b\}$ and $\mathcal{L}(w') \cap I = \emptyset.$ Therefore
\begin{displaymath}
\widetilde{t}_{sw} = \widetilde{t}_{w_{st}w'} = \widetilde{t}_{w_{st}}\widetilde{t}_{w'}  
\end{displaymath}

The hypothesis of Lemma~\ref{lemma:Z} (ii) is satisfied and we have  
\begin{displaymath}
\widetilde{t}_{w_{st}}\widetilde{t}_{w'} \in v^{-1}\mathcal{L}^s \subset v^{-1}\mathcal{L}
\end{displaymath}  
and the claim is true.

\emph{Case (iii)} In this case, $\{s, t\} = \{a, b\}$ and $\alpha, \beta \in \boldsymbol{c}$ are such that $s\boldsymbol{[s]} = \boldsymbol{[uc]} \in \textnormal{Co}(X, S)$. We can now argue as in case (ii) and the result follows.

\emph{Case (iv)} In this case, $x = sw = sw'ba$ reduced and neither $\alpha$ nor $\beta$ is in \emph{c}. We can right star reduce {\it x} to the element $sw'b$. Now, $sw'b$ and $sw'$ are both weakly complex since $\alpha, \beta \notin \boldsymbol{c}$. Thus the induction hypothesis is satisfied and we have $\widetilde{t}_{sw'b} \in v^{-1}\mathcal{L}$ and $\widetilde{t}_{sw'} \in v^{-1}\mathcal{L}$. Therefore $v\widetilde{t}_{sw'} \in \mathcal{L}$ and
\begin{displaymath}
(v\widetilde{t}_{sw'})\widetilde{t}_b = v\widetilde{t}_{sw'b} \in \mathcal{L}^b
\end{displaymath}  
by Lemma~\ref{lemma:A} (i) and Proposition 2.1.9 (ii).
By Proposition 2.1.9 (iii), $(v\widetilde{t}_{sw'b})\widetilde{t}_a \in \mathcal{L}^a$. Hence 
\begin{displaymath}
\widetilde{t}_{sw'b}\widetilde{t}_a= \widetilde{t}_{sw'ba} = \widetilde{t}_{x} \in v^{-1}\mathcal{L}^a \subset  v^{-1}\mathcal{L}
\end{displaymath}  
and the claim is true.

\emph{Case (v)} In this case, $\{s, t\} = \{b', b\}$ for $b \neq a$, and $x = sw = u'w_{st}u''a =  x'ba$ reduced. We have $ \beta \in \emph{c}$ but $\alpha \notin \emph{c}$ and $s\boldsymbol{[s]} = [\boldsymbol{u}  \boldsymbol{c} \alpha] \in \textnormal{Co}(X, S)$. Since $w_{st}$ is the longest element in the parabolic subgroup $\langle s, t \rangle$, we can use the relation $$\widetilde{t}_{w_{st}} = - \sum_{\substack{u \in \langle s, t \rangle \\ u < w_{st}}} v^{\ell(u) - m(s, t)}\widetilde{t}_u$$ to write \begin{align*} \widetilde{t}_x & = \widetilde{t}_{u'}\widetilde{t}_{w_{st}}\widetilde{t}_{u''}\widetilde{t}_a\\
& = - \widetilde{t}_{u'}\left(v^{-1}\widetilde{t}_{tw_{st}} + v^{-1}\sum_{u \leq sw_{st}} v^{\ell(u) - \ell(sw_{st})}\widetilde{t}_u\right)\widetilde{t}_{u''}\widetilde{t}_a. \end{align*}

By Lemma~\ref{lemma:M}, we have $\widetilde{t}_{u'}\widetilde{t}_u\widetilde{t}_{u''}\widetilde{t}_a \in \mathcal{L}$ where $u$ is as in the sum above. To see that $\widetilde{t}_{u'}\widetilde{t}_{tw_{st}}\widetilde{t}_{u''}\widetilde{t}_a \in \mathcal{L}$, first note that $u'tw_{st}u''$ is an expression for $x'$, which is fully commutative. If $x'$ has a reduced expression ending in $a$, since $x'b$ is weakly complex, by Theorem 3.3.8 (i) the only element in $\mathcal{R}(x')$ not commuting with $b$ is $b'$. But then $a = b'$ and $[\boldsymbol{c}a]$ is a convex chain of length $m(s, t) + 1$ and $s[\boldsymbol{s}]$ is not reduced, a contradiction. Therefore $x'$ has no reduced expression ending in $a$, $x'a$ is either fully commutative or weakly complex, and by Lemma~\ref{lemma:A} (i), $\widetilde{t}_{x'a} \in \mathcal{L}$. Putting all this together, we have $\widetilde{t}_x \in v^{-1}\mathcal{L}.$ \QED

\section{Forbidden Configurations}
The following is a collection of combinatorial results concerning inadmissible configurations in heaps of fully commutative elements. For all of \S{4.3}, let $X$ be a straight line graph with vertices $s_1, s_2, \ldots, s_n$ such that $s_i, s_j$ are adjacent if and only if $|i - j| = 1.$ 

\begin{lemma} \label{lemma:B3}
Let $X$ be a Coxeter graph of type $B_3$ and assume $m(s_1, s_2) = 4$. If $w \in W_c(X)$, then a reduced expression for $w$ has at most three occurrences of each of $s_2$ and $s_3$.
\end{lemma}
\emph{Proof.}
It follows from \cite[\S{1.8}]{Humphreys:A} that we have $w < w_0$ in the Bruhat order. This means that given a reduced expression $\boldsymbol{s}$ for $w_0$, some subexpression of $w_0$ is equal to $w$, by \cite[\S{5.10}]{Humphreys:A}. One such reduced expression for $w_0$ is $\boldsymbol{s} = s_1s_3s_2s_1s_3s_2s_1s_3s_2$. Since every expression for $w \in W_c$ has the same number of generators of each type, and $\boldsymbol{s}$ contains three occurrences of $s_2$ and three occurrences of $s_3$, the result follows. \QED

Recall from \S{3.3} and Theorem 3.3.5 that each heap $E$ is associated to a unique trace $\overline{\phi}^{-1}(E)$ , which we call the trace of $E$.

\begin{lemma} \label{lemma:M3}
Let $X$ be a straight line Coxeter graph and let $E$ be the heap of $w \in W_c(X)$.\\
(i) If $m(s_j, s_{j+1}) = 3$ for all $j \geq i$ and $[\boldsymbol{x}s_i\boldsymbol{y}s_i\boldsymbol{z}]$ is a trace for $E$, then $s_{i - 1}$ must occur in $\boldsymbol{y}$.\\
(ii) If $m(s_j, s_{j+1}) = 3$ for all $ j \leq i$ and $[\boldsymbol{x}s_i\boldsymbol{y}s_i\boldsymbol{z}]$ is a trace for $E$, then $s_{i+1}$ must occur in $\boldsymbol{y}$. 
\end{lemma}
\emph{Proof.}
Without loss of generality, we may assume that there are no occurrences of $s_i$ in $\Y$. We prove (i) by induction on $n - i$; (ii) follows similarly. If $i = n$, then we have $\overline{\phi}^{-1}(E) = [\boldsymbol{x}s_n\boldsymbol{y}s_n\boldsymbol{z}]$. Since there are no generators $s_j$ such that $j> n$, if there is no occurrence of $s_{i-1}$ in $\boldsymbol{y}$, $\boldsymbol{x}s_n\boldsymbol{y}s_n\boldsymbol{z}$ is not reduced, which completes the base case.

Assume that $i < n$. Assume by way of contradiction that there are no occurrences of $s_{i-1}$ in $\boldsymbol{y}$. Since the only other generator that does not commute with $s_i$ is $s_{i+1}$, there must be at least one occurrence of $s_{i+1}$ in $\boldsymbol{y}$, else $\boldsymbol{x}s_i\boldsymbol{y}s_i\boldsymbol{z}$ would not be reduced. If there is only one occurrence of $s_{i+1}$ in $\Y$, we have a contradiction to Theorem 3.3.8 (i) because $m(s_i, s_{i+1}) = 3$. If there is more than one occurrence of $s_{i+1}$ in $\Y$, then we are done by induction because there exist two occurrences of $s_{i+1}$ with no $s_i$ between them. \QED

\begin{lemma} \label{lemma:F5}
Let $X$ be a Coxeter graph of type $F_n$ and let $E$ be the heap of $w \in W_c(X)$. Assume that $m(s_2, s_3) = 4$.\\
(i) If $[\boldsymbol{x}s_3\Y s_3\boldsymbol{z}]$ is the trace of $E$ and $s_3, s_4$ are not in $\Y$, then $\Y$ contains a unique $s_2.$\\
The heap $E$ has no trace of the following forms:\\
(ii) $[\boldsymbol{x}(s_1s_2)\boldsymbol{y}(s_3s_2s_1s_3)\boldsymbol{z}]$ in which $s_3, s_4 \notin \boldsymbol{y}$;\\
(iii) $[\boldsymbol{x}(s_4s_3)\boldsymbol{y}(s_2s_3s_2s_4)\boldsymbol{z}]$ in which $s_4, s_5 \notin \boldsymbol{y};$\\
(iv) $[\boldsymbol{x}(s_3s_1s_2s_3)\boldsymbol{y}(s_2s_1)\boldsymbol{z}]$ in which $s_3, s_4 \notin \boldsymbol{y}$;\\
(v) $[\boldsymbol{x}(s_4s_2s_3s_2)\boldsymbol{y}(s_3s_4)\boldsymbol{z}]$ in which $s_4, s_5 \notin \Y.$
\end{lemma}
\emph{Proof.}
For part (i), $\Y$ must contain an occurrence of $s_2$ for the expression to be reduced. If $\Y$ contains two occurrences of $s_2$, this contradicts Lemma~\ref{lemma:M3} (ii).

For part (ii), we may assume $\Y$ is as short as possible, which means that only $s_1$ and $s_2$ can appear in $\Y$, since by assumption, $\Y$ has no occurrences of $s_3$ or $s_4$. If $\Y$ is empty, we have a subword of the form $s_2s_3s_2s_1s_3$ which contradicts Theorem 3.3.8 (i). If the rightmost generator of $\Y$ is $s_1$, we have $s_1s_3s_2s_1s_3$, which also contradicts Theorem 3.3.8 (i). If the rightmost generator of $\Y$ is $s_2$, we have $s_2s_3s_2s_1s_3$, which contradicts Theorem 3.3.8 (i) as well.

For part (iii), we may assume $\Y$ is as short as possible, which means that only $s_1, s_2$ and $s_3$ can appear in $\Y$, since by assumption, $\Y$ has no occurrences of $s_4$ or $s_5$. If there is an occurrence of $s_3$ in $\Y$, then by part (i), $\Y$ has no occurrences of $s_2$ to the right of the rightmost occurrence of $s_3$. But now the subword $s_3s_2s_3s_2$ contradicts Theorem 3.3.8 (i), since $m(s_2, s_3) = 4$, thus there is no occurrence of $s_3$ in $\boldsymbol{y}.$ By part (i) again, there is no occurrence of $s_2$, either. But then $\Y = s_1$, yielding another contradiction to Theorem 3.3.8 (i).

Parts (iv) and (v) are handled by symmetrical arguments to parts (ii) and (iii), respectively. The conclusion follows. \QED

\begin{lemma} \label{lemma:Fn}
Let $X$ be a Coxeter graph of type $F_n$ and let $E$ be the heap of $w \in W_c(X)$. Assume that $m(s_2, s_3) = 4$. Then $E$ cannot have a trace of the following form $$[\boldsymbol{x}(s_is_{i + 2}s_{i + 1})\boldsymbol{y}(s_{i + 1}s_{i + 2}s_i)\boldsymbol{z}]$$ in which $\boldsymbol{y}$ has no occurrences of $s_{i + 2}$ or $s_{i + 3}$.
\end{lemma}
\emph{Proof.}
The proof is by induction on $i$. Assume that $i = 1$, and that $E$ does have such a trace. Then $$[\boldsymbol{x}(s_is_{i + 2}s_{i + 1})\boldsymbol{y}(s_{i + 1}s_{i + 2}s_i)\boldsymbol{z}] = [\boldsymbol{x}(s_3s_1s_2)\boldsymbol{y}(s_2s_1s_3)\boldsymbol{z}]$$ where $s_3, s_4 \notin \boldsymbol{y}$.
We may assume $\boldsymbol{y}$ is as short as possible, which means that only $s_1$ and $s_2$ can  appear in $\boldsymbol{y}$. Now, $\boldsymbol{y}$ is non-empty, else $[\boldsymbol{x}(s_3s_1s_2)\boldsymbol{y}(s_2s_1s_3)\boldsymbol{z}]$ would not be reduced. Since $s_1s_2\boldsymbol{y}s_2s_1$ has length greater than 4 and $m(s_1, s_2) = 3$, $s_1s_2\boldsymbol{y}s_2s_1$ is not reduced, a contradiction.

Assume that $i = 2$, and that $E$ does have such a trace. Then $$[\boldsymbol{x}(s_is_{i + 2}s_{i + 1})\boldsymbol{y}(s_{i + 1}s_{i + 2}s_i)\boldsymbol{z}] = [\boldsymbol{x}(s_4s_2s_3)\boldsymbol{y}(s_3s_2s_4)\boldsymbol{z}]$$ where $ s_4, s_5 \notin \boldsymbol{y}$. We may assume that $\boldsymbol{y}$ is as short as possible, which means that only $s_1, s_2$ and $s_3$ can appear in $\boldsymbol{y}$. By Lemma~\ref{lemma:B3}, $s_2s_3\boldsymbol{y}s_3s_2$ has at most three occurrences of each of $s_2$ and $s_3.$ Suppose for a contradiction that there is an occurrence of $s_3$ in $\Y$. In this case, since there are no occurrences of $s_4$ or $s_5$ by hypothesis, there must be at least two occurrences of $s_2$ in $\boldsymbol{y}$. However, since there are now four occurrences of $s_2$ in $s_2s_3\boldsymbol{y}s_3s_2$, this contradicts Lemma~\ref{lemma:B3}. Thus it must be the case that there are no occurrences of $s_3$ in $\boldsymbol{y}$. By Lemma~\ref{lemma:F5} (i), $\Y$ contains a unique $s_2.$ But then there must be two occurrences of $s_1$ in $\boldsymbol{y}$, else the alternating product of $s_2$ and $s_3$ corresponds to a convex chain. In fact, we must have $\Y = s_1s_2s_1$. But then $s_2s_3s_1s_2s_1s_3s_2$ is not fully commutative since $m(s_1, s_2) = 3$, a contradiction.

Assume now that $i > 2$, and that $E$ does have such a trace. Note also that $m(s_i, s_{i+1}) = 3$. Then $E$ has a trace of the form $[\boldsymbol{x} s_i s_{i+2} s_{i+1}\boldsymbol{y'} s_{i+1} \boldsymbol{z'}]$
in which $\boldsymbol{y'}$ has no occurrences of $s_{i+1}$. Since $\boldsymbol{y'}$ is a subword of $\boldsymbol{y}$, $\boldsymbol{y'}$ cannot contain any occurrences of $s_{i+2}$ and $s_{i+3}$.  Since $\boldsymbol{y'}$ contains no occurrences of $s_{i+1}$, it must be the case that $\boldsymbol{y'}$ contains at least two occurrences of $s_i$, or we would contradict Theorem 3.3.8 (i).  There must also be an occurrence of $s_i$ in $\boldsymbol{z'}$, because of the subword $[s_{i+1} s_{i+2} s_i \boldsymbol{z}]$ in the original form of the trace.  We can now write the trace as
    $$[\boldsymbol{x} s_i s_{i+2} s_{i+1}\boldsymbol{y_1} s_i \boldsymbol{y_2} s_i \boldsymbol{y_3}s_{i+1}\boldsymbol{z_1} s_i \boldsymbol{z_2}],$$
where $\boldsymbol{y_1}, \boldsymbol{y_3}$ and $\boldsymbol{z_1}$ contain no occurrences of $s_i$, and $\boldsymbol{y_1}, \boldsymbol{y_2}$ and $\boldsymbol{y_3}$ contain no occurrences of $s_{i+1}, s_{i+2}$ or $s_{i+3}$. Furthermore, $s_{i+1}$ cannot appear in $\boldsymbol{z_1},$ or else we would have $[s_{i+1}\boldsymbol{z_1}s_i] = [s_{i+1}\boldsymbol{z'_1}s_{i+1}\boldsymbol{z''_1}]$ which contradicts Lemma 4.3.2 (i) applied to $[s_{i+1}\boldsymbol{z'_1}s_{i+1}]$. 

Applying the relation $[\boldsymbol{y_3}s_{i+1}] = [s_{i+1}\boldsymbol{y_3}]$, the trace is equal to $$[\boldsymbol{x} s_i s_{i+2} s_{i+1}\boldsymbol{y_1} s_i \boldsymbol{y_2}s_i s_{i+1}\boldsymbol{y_4} s_i \boldsymbol{z_2}],$$
where $\boldsymbol{y_4} = \boldsymbol{y_3 z_1}$. The generator $s_{i-1}$ must appear in $\boldsymbol{y_1}$, or the subword $[s_i s_{i+2} s_{i+1} \boldsymbol{y_1}s_i]$ would contradict Lemma 4.3.2 (i).  We can thus write $\boldsymbol{y_1} = \boldsymbol{y'_1} s_{i-1} \boldsymbol{y''_1}$,
where $s_{i-1}$ does not appear in $\boldsymbol{y''_1}$.  A similar argument applied to the subword $[s_i s_{i+1} \boldsymbol{y_4} s_i]$ shows that $s_{i-1}$ must appear in $\boldsymbol{y_4}$, and we can write $\boldsymbol{y_4} = \boldsymbol{y'_4} s_{i-1} \boldsymbol{y''_4}$, where $s_{i-1}$ does not appear in $\boldsymbol{y'_4}$.  The trace is now equal to
  $$[\boldsymbol{x} s_i s_{i+2} s_{i+1}\boldsymbol{y'_1} s_{i-1} \boldsymbol{y''_1} s_i \boldsymbol{y_2} s_i s_{i+1} \boldsymbol{y'_4} s_{i-1} \boldsymbol{y''_4} s_i \boldsymbol{z_2}].$$ Note that we have
   $$[s_{i+1} \boldsymbol{y'_1} s_{i-1} \boldsymbol{y''_1} s_i] = [\boldsymbol{y'_1} s_{i+1} s_{i-1} s_i \boldsymbol{y''_1}],$$ since $\boldsymbol{y'_1}$ has no occurrences of $s_i$ or $s_{i+2}.$ Recall that $\boldsymbol{y'_4}$ contains no occurrences of $s_{i-1}, s_{i+1}$ or $s_i$. It follows that the trace $[s_is_{i+1}\boldsymbol{y'_4} s_{i-1}]$ is equal to a trace of the form $[\boldsymbol{y'_5} s_is_{i+1}s_{i-1} \boldsymbol{y''_5}]$, where the indicated occurrences of $s_{i-1}, s_i$ and $s_{i+1}$ are identified, $\boldsymbol{y'_5}$ contains no occurrences of $s_j$ for $j > i -2$, and $\boldsymbol{y''_5}$ contains no occurrence of $s_j$ for $j < i + 2$. The inductive hypothesis can now be applied to the subword $$[\boldsymbol{y'_1} s_{i+1} s_{i-1} s_i \boldsymbol{y''_1}\boldsymbol{y_2 y'_5} s_i s_{i+1} s_{i-1} \boldsymbol{y''_5}],$$
which completes the proof.
\QED

\begin{lemma} \label{lemma:Hn}
Let $X$ be a Coxeter graph of type $H_n$ and let $E$ be the heap of $w \in W_c(X)$. Assume that $m(s_1, s_2) = 5$. Then $E$ cannot have a trace of the form $$[\boldsymbol{x}(s_is_{i + 2}s_{i + 1})\boldsymbol{y}(s_{i + 1}s_{i + 2}s_i)\boldsymbol{z}]$$ in which $\boldsymbol{y}$ has no occurrences of $s_{i + 2}$ or $s_{i + 3}$.
\end{lemma}
\emph{Proof.}
Assume that $i = 1$. Then $$[\boldsymbol{x}(s_is_{i + 2}s_{i + 1})\boldsymbol{y}(s_{i + 1}s_{i + 2}s_i)\boldsymbol{z}] = [\boldsymbol{x}(s_3s_1s_2)\boldsymbol{y}(s_2s_1s_3)\boldsymbol{z}]$$ where $s_3, s_4 \notin \boldsymbol{y}$. We may assume that $\boldsymbol{y}$ is as short as possible, which means that only $s_1$ and $ s_2$ can appear in $\boldsymbol{y}$. Now, $\boldsymbol{y}$ is non-empty, else $[\boldsymbol{x}(s_3s_1s_2)\boldsymbol{y}(s_2s_1s_3)\boldsymbol{z}]$ would not be reduced. Since $s_1s_2\boldsymbol{y}s_2s_1$ has length greater than 4 and $m(s_1, s_2) = 5$, we must have $\Y = s_1$. But this contradicts Theorem 3.3.8 (i).

Assume now that $i > 1$. We may assume that $\boldsymbol{y}$ is as short as possible, which means that only $s_j$ with $j \leq i + 1$ can appear in $\boldsymbol{y}$. Since $m(s_i, s_{i+1})=3$, we can now use the argument dealing with the case $i > 2$ from Lemma~\ref{lemma:Fn}, and the result follows. \QED

\begin{lemma}
Let $X$ be a Coxeter graph of type $\widetilde{C}_n$, ($n$ odd) and let $E$ be the heap of $w \in W_c(X)$. Then $E$ has no trace of the following forms:\\
 (i) $[\boldsymbol{x}(s_{i + 2}s_is_{i + 1})\boldsymbol{y}(s_{i + 1}s_i)\boldsymbol{z}]$ in which $\boldsymbol{y}$ has no occurrences of $s_i$ or $s_{i-1}$;\\
 (ii) $[\boldsymbol{x}(s_is_{i + 1})\boldsymbol{y}(s_{i + 1}s_is_{i + 2})\boldsymbol{z}]$ in which $\boldsymbol{y}$ has no occurrences of $s_i$ or $s_{i-1}$;\\
 (iii) $[\boldsymbol{x}(s_{i+2}s_{i + 1})\boldsymbol{y}(s_{i + 1}s_{i + 2}s_i)\boldsymbol{z}]$ in which $\boldsymbol{y}$ has no occurrences of $s_{i + 2}$ or $s_{i + 3};$\\
 (iv) $[\boldsymbol{x}(s_is_{i + 2}s_{i+1})\boldsymbol{y}(s_{i + 1}s_{i + 2})\boldsymbol{z}]$ in which $\boldsymbol{y}$ has no occurrences of $s_{i + 2}$ or $s_{i + 3}$.
\end{lemma}
\emph{Proof.}
By left and right symmetry and symmetry of the Coxeter graph, it is enough to prove (iv). The proof is by induction on $i$. Assume $i = 1$. We have $$[\boldsymbol{x}(s_is_{i + 2}s_{i + 1})\boldsymbol{y}(s_{i + 1}s_{i+2})\boldsymbol{z}] = [\boldsymbol{x}(s_3s_1s_2)\boldsymbol{y}(s_2s_3)\boldsymbol{z}]$$ where there are no occurrences of $s_3$ or $s_4$ in $\boldsymbol{y}$. We may assume $\boldsymbol{y}$ is as short as possible, which means that only $s_1$ and $s_2$ can  appear in $\boldsymbol{y}$. Now, $\boldsymbol{y}$ is non-empty, else $[\boldsymbol{x}(s_3s_1s_2)\boldsymbol{y}(s_2s_3)\boldsymbol{z}]$ would not be reduced. For the same reason, $\Y \neq s_2.$ If $\Y = s_1$, then Theorem 3.3.8 (i) is violated since $m(s_1, s_2) = 4$. Thus $s_1s_2\boldsymbol{y}s_2$ has length greater than 4, so $s_1s_2\boldsymbol{y}s_2$ is not reduced, a contradiction.

Assume now that $i > 1$.  We can write the trace in the form $$[\boldsymbol{x} s_i s_{i+2} s_{i+1} \boldsymbol{y_1} s_i \boldsymbol{y_2} s_i s_{i+1} \boldsymbol{y_4} s_i \boldsymbol{z_2}],$$
where the notation is as in the proof of Lemma~\ref{lemma:Fn}.  The proof is completed by copying
the steps of the argument used to deal with the term $\boldsymbol{y_1}$ in the case $i > 2$ in Lemma~\ref{lemma:Fn}. \QED

\section{The non-acyclic case}

 Recall from Definition 3.2.3 that $R$ is the relation defined on the vertices of a heap $E$ given by $v_1 \ R \ v_2$ if and only if  $\partial_{E}(e_0) = v_1 + v_2$ for some edge $e_0$ of $E$.
 
\begin{lemma} \label{lemma:f5}
Let $X$ be a Coxeter graph of type $\widetilde{F}_5$ and let $E$ be the heap of $w \in W_c(X)$. Let $\alpha$ be a boundary vertex in $E$. Then either $\alpha$ is an effective boundary vertex or there exists an edge $e$ in $E$ such that $\partial(e) = \alpha + \beta$, where $\beta$ is an effective boundary vertex and $\varepsilon(\alpha) \neq \varepsilon(\beta).$
\end{lemma}
\emph{Proof.}
Let the vertices of $X$ be labelled $s_0, s_1, \ldots, s_5$ and assume that $m(s_2, s_3) = 4$. By Theorem 3.4.1, there is a sequence $$\alpha_0 \ R \ \alpha_1 \ R \cdots R \ \alpha_k$$ in which $\alpha_0$ is effective. By hypothesis, we have $k > 0.$ Since $X$ is symmetric, we can deal with the case where $\varepsilon(\alpha_0) = s_2$. Since $\alpha_0$ is effective, we have $\varepsilon(\alpha_0) = s_2$ or $s_3.$ If $\varepsilon(\alpha_0) = s_0$, we must have $k = 1$ since there is no generator to the left of $s_0$, and we are done.

The only other possibility is that $\varepsilon(\alpha_1) = s_2.$ Then $E$ has a trace of the form $$[\boldsymbol{x}s_1s_2\boldsymbol{y}s_3s_2s_1s_3\boldsymbol{z}],$$ where $\Y$ has no occurrences of $s_0, s_1$ or $s_2.$ 

Suppose for a contradiction that $\Y$ contains an occurrence of $s_3.$ Then the trace is equal to $$[\boldsymbol{x}s_1s_2\boldsymbol{y_1}s_3\boldsymbol{y_2}s_3s_2s_1\boldsymbol{z}],$$ where there are no occurrences of $s_3$ or $s_2$ in $\boldsymbol{y_2}$ or $\boldsymbol{y_3}.$ Applying  Lemma~\ref{lemma:M3} (i) to $[s_3\Y_2s_3]$ gives the required contradiction.

Since $\Y$ has no occurrences of $s_3$, we have $[s_2\boldsymbol{y}] = [\boldsymbol{y}s_2]$, and  $$[\boldsymbol{x}s_1s_2\boldsymbol{y}s_3s_2s_1s_3\boldsymbol{z}] = [\boldsymbol{x}s_1\boldsymbol{y}s_2s_3s_2s_1s_3\boldsymbol{z}] = [\boldsymbol{x}s_1\boldsymbol{y}s_2s_3s_2s_3s_1\boldsymbol{z}],$$ which contradicts Theorem 3.3.8 (i).  \QED

\begin{definition} \label{def:A}
Let $E$ be a heap and let 
\begin{displaymath}
\emph{c} = \{\alpha_1, \alpha_2, \ldots, \alpha_r\} \ : \ \alpha_1 < \alpha_2 < \cdots < \alpha_r
\end{displaymath}
be a chain in {\it E}. We say \emph{c} is {\it balanced} if $\varepsilon(\alpha_1) = \varepsilon(\alpha_r)$. If \emph{c} is a balanced convex chain, we define the heap $E \backslash \emph{c}$ to be the subheap of $E$ obtained by omitting the vertices $\alpha_2, \alpha_3, \ldots, \alpha_r$. We call the heap $E \backslash \emph{c}$ the {\it contraction of E along \emph{c}}, and the number {\it r} the {\it length} of the chain.
\end{definition}

\begin{remark} \label{rem:A}
If, in the above definition, we had omitted vertices $\alpha_1, \alpha_2, \ldots, \alpha_{r-1}$ instead, we would have obtained the the same heap: the fact that \emph{c} is balanced makes the two corresponding labelled heaps isomorphic.
\end{remark}

\begin{lemma} \label{lemma:H}
Let {\it E} be a heap and let \emph{c} be a balanced convex chain of $E$.\\
(i) If the length of \emph{c} is 2, then $ \textrm{dim ker} \ \partial_{E} = \textrm{dim ker} \ \partial_{E \backslash \emph{c}} + 1$.\\
(ii) If \emph{c} is a balanced convex chain $ x < y < z$ such that $\varepsilon(x) \neq \varepsilon(y)$, then $ \textrm{dim ker} \ \partial_{E} = \textrm{dim ker} \ \partial_{E \backslash \emph{c}}$.
\end{lemma}
\emph{Proof.}
Part (i) is \cite[Lemma 2.3.4]{Green:B} and part (ii) is \cite[Lemma 2.3.5]{Green:B}. \QED

\begin{lemma} \label{lemma:J}
Let $W = W(X)$ be a star reducible Coxeter group with simply laced or complete Coxeter graph $X$ and let $\boldsymbol{s}$ be a reduced expression for $w \in W_c$. If $x = sw \notin W_c$, then the heap $ E = \overline{\phi}(s\boldsymbol{[s]})$ is acyclic. \end{lemma}
\emph{Proof.}
In the simply laced case, this is a restatement of \cite[Theorem 2.4.2 (ii)]{Green:B}, which shows that if $\overline{\phi}(s\boldsymbol{[s]})$ is not acyclic, then $\boldsymbol{[s]} = [s \UU]$ for some $\UU$, in which case $x = sw \in W_c$. If $X$ is complete, $E$ is totally ordered and $E$ will have Property P1 as long as $s[\boldsymbol{s}]$ has no trace of the form $[\boldsymbol{x}s_is_i\Y]$. This is true since $s[\boldsymbol{s}]$ is reduced. \QED

\begin{lemma} \label{lemma:4.4.6}
Let $W(X)$ be a star reducible Coxeter group, let $\boldsymbol{s}$ be a reduced expression for $w \in W_c$ and let $x = sw \notin W_c$ for some $s \in S$. Let $E$ be the heap of $w$ and assume $\overline{\phi}(s\boldsymbol{[s]})$ is not acyclic. If a boundary vertex $\alpha$ occurs in a sequence $\alpha_0 \ R \ \alpha_1 \ R \ \cdots \ R \ \alpha_k = \alpha$ where $\alpha_0, \alpha_1, \ldots, \alpha_k = \alpha$ are  boundary vertices and $\alpha_0$ is an effective boundary vertex, then $\varepsilon(\alpha_i) \neq \varepsilon(\alpha_{i + 1})$, except possibly if $i = k - 1.$
\end{lemma}
\emph{Proof.} The simply laced and complete Coxeter graphs are dealt with in Lemma~\ref{lemma:J}. Thus $X$ is a Coxeter graph of type $B_n,\ F_n, \ \widetilde{C}_n, (n \ \textnormal{odd}), \ H_n$, or $\widetilde{F}_5.$ This means that if $e_1 = (x_1, y_1)$ and $e_2  = (x_2, y_2)$ are edges in $E$ with $\partial(e_1) = \alpha_l + \alpha_{l+1}$ and $\partial(e_2) = \alpha_{l+1} + \alpha_{l+2}$, then the subscripts of the Coxeter generators $\varepsilon(x_1), \varepsilon(x_2)$ differ by 0 or 2.

Suppose we are in type $\widetilde{C}_n, (n \ \textnormal{odd})$, and let $l$ be minimal such that $\varepsilon(\alpha_l) = \varepsilon(\alpha_{l+1})$ as in the statement. Let $e$ be an edge of $E$ such that $\partial(e) = \alpha_l + \alpha_{l+1}$. By symmetry, we may assume that $e = (x, y)$ with $\varepsilon(x) = \varepsilon(y) = s_{j-1}$. This cannot happen without contradicting Lemma 4.3.6 (i) or (iii) with $i = j-1$. This completes the proof in type $\widetilde{C}_n, (n \ \textnormal{odd})$, which contains type $B_n$ as a subcase.

The case where $X$ is a Coxeter graph of type $\widetilde{F}_5$ is dealt with in Lemma~\ref{lemma:f5}, which shows that $\varepsilon(\alpha_i) = \varepsilon(\alpha_{i+1})$ is impossible. 
Thus the present argument reduces to types $F_n$ and $H_n$.

\emph{Case (i).} Let $X$ be a Coxeter graph of type $F_n$. Suppose for a contradiction that we have $\alpha_l \ R \ \alpha_{l+1} \ R \ \alpha_{l+2}$ with $\varepsilon(\alpha_l) = \varepsilon(\alpha_{l+1}) = \varepsilon(\alpha_{l+2}) = s_j$. Then there are edges  $e_1 = (x_1, y_1)$ and $e_2  = (x_2, y_2)$ in $E$ with $\partial(e_1) = \alpha_l + \alpha_{l+1}$ and $\partial(e_2) = \alpha_{l+1} + \alpha_{l+2}$. We must have $\varepsilon(x_1) = s_{j-1}$ and $\varepsilon(x_2) = s_{j+1}$ or vice versa. We deal with the first case; the other case is symmetric. The trace of $E$ is now of the form $$[\boldsymbol{x}s_{j+1}s_j\boldsymbol{y_1}s_{j-1}s_js_{j+1}\boldsymbol{y_2}s_js_{j-1}\boldsymbol{z}]$$ or $$[\boldsymbol{z}s_{j-1}s_j\boldsymbol{y_2}s_{j+1}s_js_{j-1}\boldsymbol{y_1}s_js_{j+1}\boldsymbol{x}],$$ in which $\boldsymbol{y_1}$ has no occurrences of $s_j$ or $s_{j+1}$ and $\boldsymbol{y_2}$ has no occurrences of $s_j$ or $s_{j-1}.$

Let $t$ be the unique index with $m(s_t, s_{t+1}) > 3.$ If $t \geq j$, then the subword $[s_j\boldsymbol{y_1}s_{j-1}s_j]$ or $[s_js_{j-1}\boldsymbol{y_1}s_j]$ contains no occurrences of $s_{j+1}$, which contradicts Lemma 4.3.2 (ii). If $t < j$, then the subword $[s_js_{j+1}\boldsymbol{y_2}s_j]$ or $[s_j\boldsymbol{y_2}s_{j+1}s_j]$ contradicts Lemma 4.3.2 (i). 

Now suppose $\varepsilon(\alpha_1) = \varepsilon(\alpha_0) = s_j$. If $j = 2$, this contradicts Lemma 4.3.3 (ii) or (iv); if $j=3$, this contradicts Lemma 4.3.3 (iii) or (v). It is now enough to show we cannot have $s_j = \varepsilon(\alpha_l) = \varepsilon(\alpha_{l+1}) \neq \varepsilon(\alpha_{l +2}) = \varepsilon(\alpha_{l-1}) = s_{j'}$.

Assume first that $j' = j+2.$ In this case the trace is of the form $$[\boldsymbol{x}s_{j-1}s_{j+1}s_js_{j+2}s_{j+1}\Y s_{j+1}s_js_{j+2}s_{j-1}s_{j+1}\boldsymbol{z}]$$ where $\Y$ has no occurrences of $s_{j-1}$ or $s_j$. Now the subword $s_js_{j+2}s_{j+1}\Y s_{j+1}s_j$ contradicts Lemma 4.3.2 (i) unless $j=2$, in which case we contradict Lemma 4.3.3 (ii) or (iv).

We now deal with the case $j' = j - 2$. There exist distinct edges $e_1 = (x_1, y_1)$,  $e_2  = (x_2, y_2)$, and $e_3 = (x_3, y_3)$ in $E$ with $\partial(e_1) = \alpha_{l-1} + \alpha_l$, $\partial(e_2) = \alpha_{l+1} + \alpha_l$, and $\partial(e_3) = \alpha_{l+1} + \alpha_{l+2}$. Assume without loss of generality that $x_3 < y_1$. The formula for $\partial(e_3)$ shows that the trace of $E$ has subwords of the form $s_{j-2}s_js_{j-1}$ and $s_{j-1}s_js_{j-2}$, where the first subword occurs to the left of the second subword. If the occurrences of $s_{j-1}$ are the same, the trace has a subword of the form $$[s_{j-2}s_js_{j-1}s_js_{j-2}] = [s_js_{j-2}s_{j-1}s_{j-2}s_j],$$ which is impossible as either $m(s_{j-2}, s_{j-1}) = 3$ or $m(s_j, s_{j-1}) = 3$. If the occurrences are distinct, then the trace is of the form $$[\boldsymbol{x}s_{j-1}s_{j+1}s_{j-2}s_js_{j-1}\Y s_{j-1}s_js_{j-2}s_{j+1}s_{j-1}\boldsymbol{z}]$$ where $\Y$ has no occurrences of $s_j$ or $s_{j+1}$, which contradicts Lemma~\ref{lemma:Fn} with $i = j -2$. This completes the proof in type $F_n$.\\
\emph{Case (ii).} Let $X$ be a Coxeter graph of type $H_n$ and let $s_j = \varepsilon(\alpha_0)$; note that $j = 1$ or $2$. If $\varepsilon(\alpha_1) = \varepsilon(\alpha_0) = s_j$, then we cannot have $j=1$ because $s_1$ is the leftmost vertex. If $j=2$, then we must have $\partial(e)= \alpha_0 + \alpha_1$ for some edge $e = (x, y)$ with $\varepsilon(x) = s_3.$ If $k > 1$, we must have $\varepsilon(\alpha_2) = s_1$ (as opposed to $s_3$). This implies that $\varepsilon(\alpha_0) = \varepsilon(\alpha_1) = \varepsilon(\alpha_2)$, which was proved impossible above. 

We may therefore assume that $\varepsilon(\alpha_0) \neq \varepsilon(\alpha_1).$ It is now enough to show that we cannot have $$s_j = \varepsilon(\alpha_l) = \varepsilon(\alpha_{l+1}) \neq \varepsilon(\alpha_{l +2}) = \varepsilon(\alpha_{l-1}) = s_{j'}.$$ We now argue as in Case (i) above, except at the final step, where we invoke Lemma~\ref{lemma:Hn} instead of Lemma~\ref{lemma:Fn}. (In this case, the situation where $j=2$ in Lemma 4.3.3 cannot arise.)\QED

 \begin{lemma} \label{lemma:F}
Let $W(X)$ be a star reducible Coxeter group, let $\boldsymbol{s}$ be a reduced expression for $w \in W_c$ and let $x = sw \notin W_c$ for some $s \in S$. Let $E$ be the heap of $w$ and assume $\overline{\phi}(s\boldsymbol{[s]})$ is not acyclic. Let $\boldsymbol{c} = s^{(0)} < t^{(1)} < s^{(1)} < t^{(2)} < \cdots < s^{(k)}$ (or $t^{(k)}$, depending on the parity of $m(s, t)$) be the convex chain in $s\boldsymbol{[s]}$ described in Theorem 3.3.8 (i), where $s\boldsymbol{[s]} = \boldsymbol{[x'cx'']}$ and $\boldsymbol{c'} = \boldsymbol{c} \backslash \{\alpha_r\}$. Then the heap $\overline{\phi}(\boldsymbol{x'c'x''})$ is acyclic.
\end{lemma}
\emph{Proof.}
If $X$ is simply laced or complete, then the heap $\overline{\phi}(s\boldsymbol{[s]})$ is acyclic by Lemma~\ref{lemma:J}, so we may assume $X$ is a straight line graph. Let $E_0 = \overline{\phi}(s\boldsymbol{[s]})$ and let $s^{(0)}$ denote the new vertex in $E_0$. Since $E$ is acyclic by Lemma 3.3.18 (i) (iv), Lemma 4.1.1 shows that $\textrm{dim ker} \ \partial_{E_0} = 1.$ Now $e_0= (s^{(0)}, s^{(1)})$ is an edge such that $\partial(e_0) = t^{(1)}$. Since $E_0$ is not acyclic, there exist edges $e_0, e_1, \ldots e_r$ in $E_0$ and scalars $\lambda_0, \lambda_1, \ldots, \lambda_r$, not all equal to 0, such that 
\begin{displaymath}
\partial_{E_0} \left( \sum_{i=0}^r\lambda_i e_i \right) = 0,
\end{displaymath}
where $e_0=  (s^{(0)}, s^{(1)})$, so $\partial_E(e_0) = t^{(1)}$. Since $E$ is acyclic and $E$ is a convex subheap of $E_0$ by Lemma 3.3.9, Lemma 3.3.12 shows that $\lambda_0 \neq 0$, and without loss of generality we may assume $\lambda_0 = 1.$ Now we have that
\begin{displaymath}
\partial_{E_0} \left(-\sum_{i=1}^r\lambda_i e_i \right) = t^{(1)} 
\end{displaymath}
and similarly,
\begin{displaymath}
\partial_E \left(-\sum_{i=1}^r\lambda_i e_i \right) = t^{(1)} 
\end{displaymath}
by Lemma 3.3.12. Therefore $t^{(1)}$ is a boundary vertex in $E$.

By Theorem 3.4.1, there is a chain of linear equivalences $\alpha_0 \ R \ \alpha_1 \ R \cdots \ R \ \alpha_k$, where $\alpha_k = t^{(1)}$, $\alpha_0, \alpha_1, \ldots, \alpha_k$ are boundary vertices in $E$, and $\alpha_0$ is effective. By Lemma 4.4.6, the set $L = \{ \varepsilon(\alpha_0), \varepsilon(\alpha_1), \ldots, \varepsilon(\alpha_k)\}$ has cardinality $k+1$ or $k$, and in the latter case we have
$\varepsilon(\alpha_{k-1})= \varepsilon(\alpha_k)$. The proof of Lemma 4.4.6 shows that the latter situation can only occur in type $F_n$ or $H_n$.

\emph{Case (i).} Suppose first that $|L| = k+1$. In this case, $\alpha_0, \alpha_1, \ldots, \alpha_k$ correspond to the labels $s_2, s_4, \ldots, s_{k+2}$ in Figure 8 below; this implies that $m(s, t) = 3.$
\vspace{0.5 in}
\begin{center}

\includegraphics{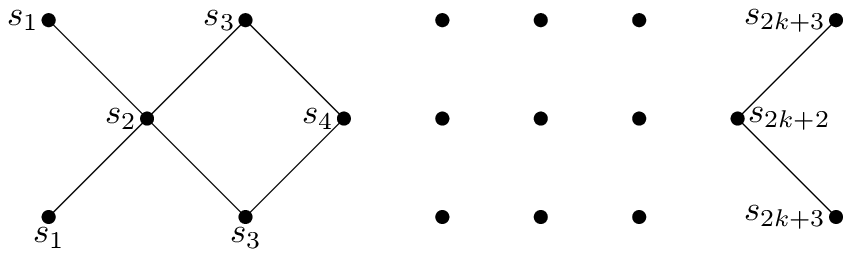}

Figure 8
\end{center}

By Lemma~\ref{lemma:H}, if we contract the chain $\emph{c} = s_1 < s_2 < s_1$, then $ \textrm{dim ker} \ \partial_{E_0} = \textrm{dim ker} \ \partial_{E_0 \backslash \emph{c}} = 1$. If we repeat this process by contracting the chain $\boldsymbol{c'} = s_3 < s_4 < s_3$ in $E_0 \backslash \emph{c}$, then $ \textrm{dim ker} \ \partial_{E_0 \backslash \emph{c}} = \textrm{dim ker} \ \partial_{{(E_0 \backslash \emph{c})} \backslash \boldsymbol{c'}} = 1$.
By iterating this process, we wind up with the heap $E'$ in Figure 9 below.
\vspace{.5 in}
\begin{center}

\includegraphics{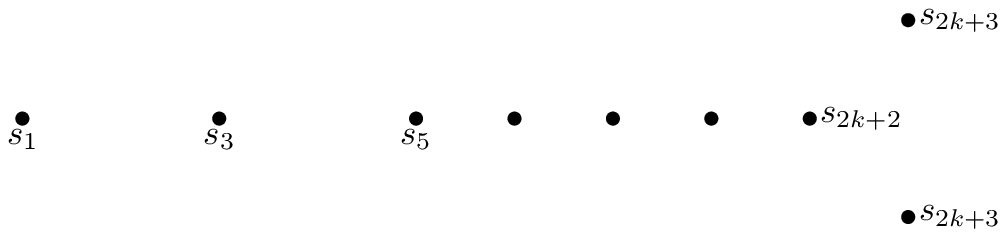}
Figure 9
\end{center}

Contracting along $\boldsymbol{c''} = s_{2k+3}  < s_{2k+3}$ then gives $$ \textrm{dim ker} \ \partial_{E'_0 \backslash \boldsymbol{c''}} = \textrm{dim ker} \ \partial_{E'_0} - 1.$$ Since $\textrm{dim ker} \ \partial_{E_0} = 1$, $E'_0  \backslash \boldsymbol{c''}$ is acyclic. Since the lower occurrence, $\gamma = s^{(1)}$, of $s_{2k+3}$ is only involved in the last contraction, $\textrm{dim ker} \ \partial_{{E_0} \backslash s^{(1)}} < \textrm{dim ker} \ \partial_{E_0} = 1.$ Thus $E_0 \backslash s^{(1)}$ is acyclic, as required.

\emph{Case (ii).}
Now assume $|L| = k$ and $X$ is of type $H_n$ or $F_n$. Assume first that the chain $\alpha_0 \ R \ \alpha_1 \ R \ \cdots \ R \ \alpha_k = \alpha$ and its corresponding edges is as shown up to symmetry in Figure 10 below, where the four occurrences of $s_{k-1}$ are distinct and $\varepsilon(\alpha_l) = s_{2l+2}$ for $l < k$; this implies that $m(s, t) = 3$. We proceed as in case (i) starting with the chain $s_1 < s_2 < s_1$, removing the $\alpha_i$ in order starting with $\alpha_0$ and ending with the chain $s^{(0)}  < s^{(1)}$. Since $s^{(1)}$ is not involved with any earlier contractions, we argue as in case (i) that $E_0 \backslash s^{(1)}$ is acyclic.
\vspace{0.5 in}
\begin{center}

\includegraphics{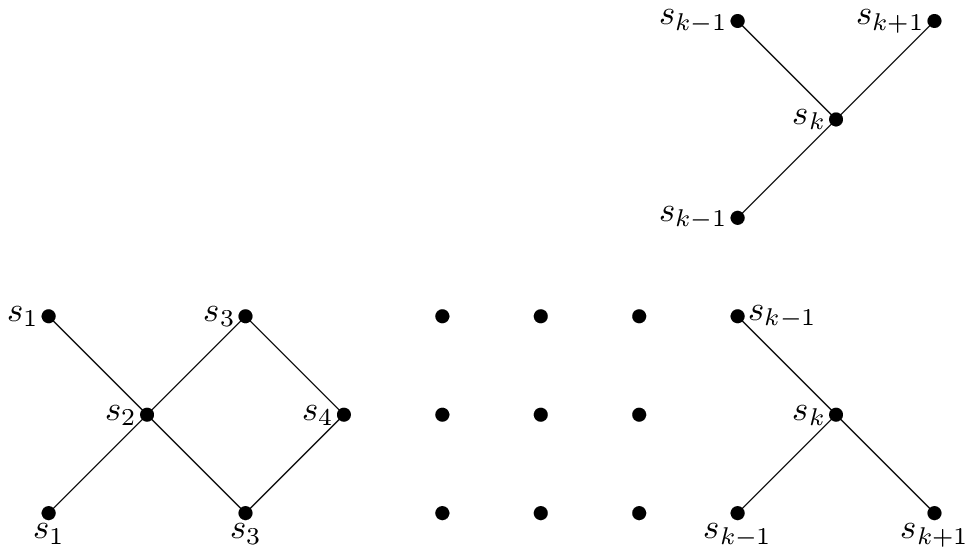}

Figure 10
\end{center}

\emph{Case (iii).}
The other possibility is that $|L| = k$ and $X$ is of type $H_n$ or $F_n$, but the occurrences of $s_{k-1}$ in Figure 10 are the same. This situation is shown in Figure 11 below. 
\vspace{0.5 in}
\begin{center}

\includegraphics{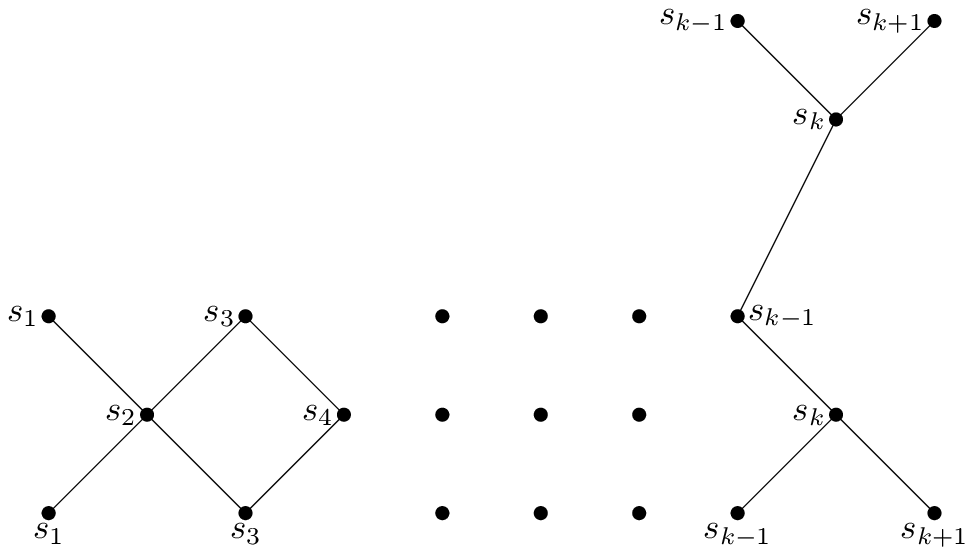}

Figure 11
\end{center}

In this case, the occurrences of $s_k$ shown are $t^{(1)}$ and $t^{(2)}$, which implies that $m(s_{k-1}, s_k) > 3$. Since $\alpha_0$ is effective, $m(s_1, s_2) > 3$. The only way this can happen is if $k=2$ and $X$ is of type $H_n$. This situation is shown in Figure 12, where the bottom left $s_1$ is $s^{(2)}$ and $m(s, t)=5.$ 
\vspace{0.5 in}
\begin{center}

\includegraphics{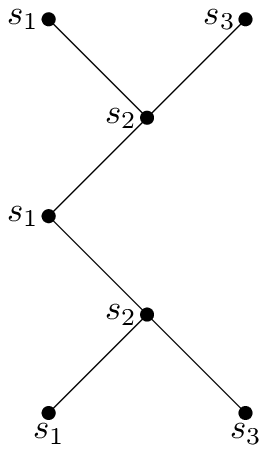}

Figure 12
\end{center}

In this case, we perform a contraction along $\boldsymbol{c}= s^{(0)} < t^{(1)} < s^{(1)}$, followed by a contraction along $\boldsymbol{c'} = s_3 < s_2 < s_3$. Since $s^{(2)}$ is not involved in either contraction, we can argue as before to show that $\textrm{dim ker} \ \partial_{{E_0} \backslash s^{(2)}} = 0,$ which completes the proof. \QED

\begin{theorem} \label{theorem:G} (\emph{Property W})
Let $W(X)$ be a star reducible Coxeter group with Coxeter graph {\it X}. Let $w \in W_c$ and let $x = sw \notin W_c$ for some $s \in S$. Then $\widetilde{t}_x \in v^{-1}\mathcal{L}$.
\end{theorem}
\emph{Proof.}\\
%"Let E be the heap of w and let [s] be the trace of E.  If \bar\phi(s[s])
%is acyclic, the conclusion follows by Lemma 4.2.1.  In particular, if X
%is simply laced or complete, we are done by Lemma 4.4.5, so we may assume
%from now on that \bar\phi(s[s]) is not acyclic."
Let $E$ be the heap of $w$ and let $[\boldsymbol{s}]$ be the trace of $E$. If $\overline{\phi}(s\boldsymbol{[s]})$ is acyclic, the conclusion follows by Lemma 4.2.1. In particular, if $X$ is simply laced or complete we are done by Lemma 4.4.5,  so we may assume from now on that $\overline{\phi}(s\boldsymbol{[s]})$ is not acyclic.

Assume that $x = sw = u'w_{st}u''$ reduced, where $w_{st}$ is the longest element in the parabolic subgroup $\langle s, t \rangle$. We can use the relation $$\widetilde{t}_{w_{st}} = - \sum_{\substack{u \in \langle s, t \rangle \\ u < w_{st}}} v^{\ell(u) - m(s, t)}\widetilde{t}_u$$ to write \begin{align*} \widetilde{t}_x & = \widetilde{t}_{u'}\widetilde{t}_{w_{st}}\widetilde{t}_{u''}\\
& = - \widetilde{t}_{u'}\left(v^{-1}\widetilde{t}_{tw_{st}} + v^{-1}\sum_{u \leq sw_{st}} v^{\ell(u) - \ell(sw_{st})}\widetilde{t}_u\right)\widetilde{t}_{u''}. \end{align*}

Arguing as in Case (v) of the proof of Lemma~\ref{lemma:C}, we find that all of the terms in the expression for $\widetilde{t}_x$ arising from $u \leq sw_{st}$ lie in $v^{-1}\mathcal{L}$. Now write $s[\boldsymbol{s}]=\boldsymbol{x'cx''}$ and $\boldsymbol{x'c'x''}$ as in Lemma~\ref{lemma:F}. By Lemma~\ref{lemma:F}, $\overline{\phi}([\boldsymbol{x'c'x''}])$ is acyclic. By Lemma 4.1.6, $\widetilde{t}_{\boldsymbol{x'c'x''}} \in \mathcal{L}$. Now, $$\widetilde{t}_{\boldsymbol{x'c'x''}}= \widetilde{t}_{u'}\widetilde{t}_{tw_{st}}\widetilde{t}_{u''},$$ so $v^{-1}\widetilde{t}_u\widetilde{t}_{tw_{st}}\widetilde{t}_{u''} \in v^{-1}\mathcal{L},$ which completes the proof. \QED

\bibliographystyle{plain} 
\bibliography{Articles} 
\end{document}